\documentclass[10pt]{article}
\usepackage[dvips]{graphics}
\usepackage{multirow}
\usepackage{array}
\usepackage{epsfig,rotating}
\usepackage{amssymb,amsmath, mathabx}
\usepackage{latexsym}
\usepackage[usenames]{color}

\numberwithin{equation}{section}

\newcommand{\bd}{\mathbf}

\newcommand{\bZ}{\mathbf{Z}}
\newcommand{\bz}{\mathbf{z}}
\newcommand{\bx}{\mathbf{x}}
\newcommand{\bY}{\mathbf{Y}}
\newcommand{\bS}{\mathbf{S}}
\newcommand{\bT}{\mathbf{T}}

\def\Var{\mathrm{Var}}
\newtheorem{theorem}{Theorem}[section]
\newtheorem{lemma}{Lemma}[section]

\newcommand\eop{\hfill$\blacksquare$}

\begin{document}

\title{Largest Eigenvalues of  Principal Minors of Deformed Gaussian Orthogonal Ensembles and Wishart Matrices}

\author{Tiefeng Jiang$^{1}$
and  Yongcheng Qi$^2$\\
The Chinese University of Hong Kong at Shenzhen and
the University of Minnesota Duluth}

\date{}
\maketitle

\footnotetext[1]{School of Data Science, the Chinese University of
Hong Kong, Shenzhen, Guangdong, 518172, P.R. China.
The research of Tiefeng Jiang was supported in part by UDF01003320 and UF02003320.}

\footnotetext[2]{Department of Mathematics and Statistics,
University of Minnesota Duluth, MN 55812, USA.
The research of Yongcheng Qi was supported in part by NSF Grant
DMS-1916014.}

\begin{abstract}


\noindent Consider a high-dimensional Wishart matrix $\bd{W}=\bd{X}^T\bd{X}$
 where the entries of $\bd{X}$ are i.i.d. random variables with mean zero, variance one, and  a finite fourth moment $\eta$. Motivated by problems in signal processing and high-dimensional statistics, we study the maximum of the largest eigenvalues of any two-by-two principal minors of $\bd{W}$. Under certain restrictions on the sample size and the population dimension of $\bd{W}$, we obtain the limiting distribution of the maximum, which follows the Gumbel distribution when $\eta$ is between 0 and 3, and a new distribution when $\eta$ exceeds 3. To derive this result, we first address a simpler problem on a new object named a deformed Gaussian orthogonal ensemble (GOE). The Wishart case is then resolved using results from the deformed GOE and a high-dimensional central limit theorem. Our proof strategy combines the Stein-Poisson approximation method, conditioning, U-statistics, and the Hájek projection. This method may also be applicable to other extreme-value problems. Some open questions are posed.

\end{abstract}



\newpage

\section{Introduction}


Driven by applications in high-dimensional statistics and signal processing, Cai {\it et al.}~\cite{CJL2021} investigated the extreme eigenvalues of the principal minors of a Wishart matrix. The details are described next.

Let $\{x_{ij}; i\geq 1, j\geq 1\}$ be i.i.d. random variables with mean $0$ and variance $1$.
Write $\bd{X}=(x_{ij})_{n\times p}=(\bx_1, \cdots,\bx_p).$ Then $\bx_1, \cdots,  \bx_p$ are column vectors of $\bd{X}.$ The $p\times p$ matrix $\bd{W}=(w_{ij})_{1\le i,j\le
p}:=\bd{X}^T\bd{X}$ is referred to as a Wishart matrix in literature. For $m\geq 1$, let $S =\{i_1, \cdots, i_m\}\subset \{1,2,\cdots, p\}$ with the size of $S$ being $m$ and an $n\times m$ matrix $\bd{X}_S=(\bx_{i_1},\cdots, \bx_{i_m})$. Then $\bd{W}_S=\bd{X}_S^T\bd{X}_S$ is a $m\times m$ principal minor of $\bd{W}.$
 Denote by $\lambda_1(\bd{W}_S)\geq \cdots \geq \lambda_m(\bd{W}_S)$ the eigenvalues of $\bd{W}_S$. The quantity of interest is the maximum of the largest eigenvalues of all  $m\times m$ principal minors of $\bd{W}$. Specifically,  we will work on
\begin{equation}\label{NEWS}
T_{n,p, m}:= \max_{S\subset\{1, \cdots, p\},\, |S|=m}\lambda_1(\bd{W}_S)
\end{equation}
where $|S|$ denotes the cardinality of the set $S$.

A simpler but useful analogue of \eqref{NEWS} is that the Wishart matrix $\bd{W}$ is replaced by a so-called deformed Gaussian orthogonal ensemble (GOE). Let us explain this setting. Notice $w_{ij}=\sum_{k=1}^nx_{ki}x_{kj}$ for any $i,j$ and in particular, $w_{ii}=\sum_{k=1}^nx_{ki}^2$. As a sum of i.i.d. random variables, we have from the classical  central limit theorem that
\begin{equation}\label{chifan}
\frac{1}{\sqrt{n}}\big(w_{ij}-n\delta_{ij}\big) \Rightarrow
\begin{cases}
N(0, 1), &\text{if $i\ne j$};\\
N(0, \xi), &\text{if $i=j$}
\end{cases}
\end{equation}
weakly, where $\xi=\mbox{Var}(x_{11}^2)$,  $\delta_{ij}=1$ if $i=j$ and $0$ otherwise. Trivially, the left hand side of \eqref{chifan} is uncorrelated to
$(1/\sqrt{n})\big(w_{kl}-n\delta_{kl}\big)$ if $\{k, l\} \ne\{i, j\}.$  Based on this, we propose a new  random matrix model $\bZ=(\bz_{ij})_{1\le i,j\le p}$. It is a symmetric matrix,  $\{\bz_{ij}, 1\le i\leq j\le p\}$ are
independent normal random variables with mean zero,
$\Var(\bz_{ij})=1$ for different $i,j$ and $\Var(\bz_{ii})=\xi \geq 0$ for each $i.$ Of course, the matrix $\bZ$ is the Gaussian
orthogonal ensemble (GOE) when $\xi=2$ in the literature of random matrix theory. For this reason, we call $\bZ$ a deformed GOE associated with a parameter $\xi\geq 0.$
The analogue of \eqref{NEWS} becomes
\begin{equation}\label{NEWS1}
\tilde{T}_{p, m}= \max_{S\subset\{1, \cdots, p\},\, |S|=m}\lambda_1(\bd{Z}_S)
\end{equation}
where $\bd{Z}_S$ is the $m\times m$ principal minor of $\bd{Z}$ corresponding to $S =\{i_1, \cdots, i_m\}\subset \{1,2,\cdots, p\}$. Superficially, as observed above, $T_{n,p, m}$ and $\tilde{T}_{p, m}$ are close to each other. In fact, in this paper we will first work on  $\tilde{T}_{p, m}$ then we get the same conclusion on $T_{n,p, m}$ by a high-dimensional central limit theorem.

Under a Gaussian assumption on $x_{11}$,
Cai {\it et al.}~\cite{CJL2021}
obtained that
\begin{equation}\label{CJL1}
\frac{T_{n,p, m}-n}{\sqrt{n}}-2\sqrt{m\ln p}\to 0 ~~\mbox{ in
probability}
\end{equation}
as $n,p$ go to infinity for any fixed $m\geq 1$. Furthermore, for the GOE ($\xi=2$) they also proved that
\begin{eqnarray}\label{callback}
\tilde{T}_{p, m} -2\sqrt{m\log p} \to 0~~\mbox{ in
probability}
\end{eqnarray}
as $p\to\infty$.  A few questions were asked and discussed in Remark 4 from Cai {\it et al.}~\cite{CJL2021}. For example, (a) what are the limiting distributions of $T_{n,p, m}$ and  $\tilde{T}_{p, m}$? (b) Whether the Gaussian assumption on $x_{ij}$ required in \eqref{CJL1} can be removed? They also conjecture that $(T_{n,p, m}-n)/\sqrt{n}$ and $\tilde{T}_{p, m}$ have similar behavior if $\xi=\mbox{Var}(x_{11}^2) \leq 2.$ In this paper we will study the problems for every $\xi\geq 0$.


Recently, Feng {\it et al.} \cite{Feng2022}  obtained a result in italics next.
{\it For GOE ($\xi=2$), we have the weak convergence
\begin{equation}\label{cozy}
\tilde{T}_{p, m}^2-4m \ln p - 2(m-2) \ln\ln p \overset{d}{\longrightarrow} Y
\end{equation}
as $p\to\infty$ for fixed $m$, where the random variable $Y$ has the Gumbel distribution function
\begin{eqnarray*}
F_Y(y)=\exp\big(-c_me^{-y/4}\big), ~~~ y \in \mathbb{R}.
\end{eqnarray*}
Here, the constant
\begin{eqnarray*}
c_m=\frac{(2m)^{(m-2)/2}K_m}{(m-1)!2^{3/2}\Gamma(1+m/2)}
\end{eqnarray*}
where $K_m=\mu(\mathcal{S}_m)$ is the probability of the event
\begin{eqnarray*}
 \mathcal{S}_m:=\Big\{x\in S^{m-1}:\, \sum_{j\in \beta}x_j^2\leq \sqrt{k/m},\, \forall \beta \subset \{1, \cdots, m\}~\mbox{with}~ \forall 1\leq |\beta|=k<m\Big\}
\end{eqnarray*}
under the uniform measure $\mu$ on the unit sphere $S^{m-1}.$ In particular,
\begin{eqnarray*}
 c_1=\frac{1}{2\sqrt{\pi}}, ~~~ c_2=-\frac{1}{2\sqrt{2}} + \frac{\sqrt{2}}{\pi}\arcsin \big(2^{-1/4}\big).
\end{eqnarray*}
}
\noindent The authors derived \eqref{cozy} by using the Poisson approximation method as well as the explicit density function of the eigenvalues of the GOE ($\xi=2$). Comparing their result \eqref{cozy} to questions (a) and (b) above, in this paper we will first investigate $\tilde{T}_{p, m}$ with $m=2$ and $\xi=\Var(\bz_{11})$ being an arbitrary non-negative number. By this step and a high-dimensional central limit theorem, we obtain a theorem  on  $T_{n,p, m}$ with $m=2$ without the assumption that $x_{11}$ is Gaussian. We discover that the limiting distributions of both $\tilde{T}_{p, m}$ and $T_{n,p, m}$ are Gumbel distributions for $\xi \in [0, 2]$. However, for $\xi>2$, their limiting distributions are no longer the Gumbel distribution. Instead they are a different class of distributions with an explicit formula; see \eqref{G-xi}. To our knowledge, it is a new probability distribution.  A further understanding of the distribution needs extra effort. It is interesting to see that our Theorem \ref{main1} for GOE ($\xi=2$) agree with \eqref{cozy} as $m=2$.

Now we comment on the connection and difference between our results in this paper and \eqref{cozy} from \cite{Feng2022}. Clearly, \cite{Feng2022} studies the behavior of $\tilde{T}_{p, m}$ as $p\to\infty$ with $\xi=2$ and any $m\geq 1$. In this paper we first explore $\tilde{T}_{p, m}$ as $p\to\infty$ with arbitrary $\xi\geq 0$ and  $m=2$. Then we get the same conclusion on our initial target $T_{n,p, m}$ as $m=2$ under certain conditions on $n$ and $p$. The authors from \cite{Feng2022} did not study $T_{n,p, m}$. It is exciting to see the results from  \cite{Feng2022} on GOE and ours on the deformed GOE complement each other although the two methods of deriving them are very different. In particular, there is no explicit formula for the joint density function of the eigenvalues of a deformed GOE for $\xi \ne 2.$ As a result we have no way to manipulate the eigenvalues via its density function as done in \cite{Feng2022}.

The organization of the remaining paper is as follows. In section \ref{mainresults} we state a result on the limiting distributions of $\tilde{T}_{p, m}$ in \eqref{NEWS1} on a deformed GOE and  that of $T_{n,p, m}$ in \eqref{NEWS} relevant to Wishart matrices. The method of our proofs and open problems are also discussed in this section. The proofs of main results and discussions are presented  in  Section \ref{proof:main}. The proofs of all technical lemmas are placed at Section \ref{lemmas}.

\section{Main results}\label{mainresults}




Consider a $p\times p$ deformed Gaussian orthogonal ensemble $\bZ=(\bz_{ij})_{1\le i,j\le p}$ defined in \eqref{NEWS1}. That is, $\bZ$ is  a symmetric random matrix in which $\{\bz_{ij}, 1\le i\leq j\le p\}$ are
independent normal random variables with mean zero and
$\Var(\bz_{ij})=1$ for $1\le i<j\le p$ and $\Var(\bz_{ii})=\xi \geq 0$
for $1\le i\le p$. Of course, the case $\xi=0$ corresponds to that all diagonal entries $\bz_{ii}$ are equal to $0$.  The matrix $\bZ$ reduces to the Gaussian
Orthogonal Ensemble (GOE) when $\xi=2$ as aforementioned.

Denote $l_{ij}$ as the largest eigenvalue of the $2\times 2$ principal minor
\[
\bZ_{\{i,j\}}=
\left(
          \begin{array}{cc}
            \bz_{ii} & \bz_{ij} \\
            \bz_{ij} & \bz_{jj} \\
          \end{array}
        \right)
\]
for all $1\le i< j\le p$, that is, $l_{ij}=\lambda_1(\bZ_{\{i,j\}})$. It is easily seen that
\begin{equation}\label{Dimsum}
l_{ij}=\sqrt{\bz_{ij}^2+\frac{1}{4}(\bz_{ii}-\bz_{jj})^2}+
\frac1{2}(\bz_{ii}+\bz_{jj})
\end{equation}
for all $1\le i< j\le p$. Obviously, $\tilde{T}_{p, m}=\max_{1\le i<j\le
p}l_{ij}$ by \eqref{NEWS1}.


Set $\alpha_p=\sqrt{2\ln p}$ and $\beta_p=\alpha_p-\alpha_p^{-1}\ln(\sqrt{2\pi}\alpha_p)$. Define
 \[
A_{p,\xi}=\left\{
  \begin{array}{ll}
    \frac{2}{\sqrt{2+\xi}}\alpha_p, & \hbox{ if }~ 0\leq \xi<2; \\
    \alpha_p, & \hbox{ if }~ \xi=2; \\
    \frac{2+\sqrt{\xi}}{\xi+\sqrt{\xi}}\alpha_p, & \hbox{ if }~ \xi>2
  \end{array}
\right.
\]
and
\[
B_{p,\xi}=\left\{
  \begin{array}{ll}
    \sqrt{2+\xi}\alpha_p-\alpha^{-1}_p\frac{\sqrt{2+\xi}}{2}\ln\big(\sqrt{2\pi(2-\xi)}\alpha_p\big), & \hbox{ if } ~0\leq \xi<2; \\
    2\alpha_p-\alpha^{-1}_p\ln\frac{\sqrt{2}\pi}{\arcsin(\sqrt{2}-1)}, & \hbox{ if }~\xi=2; \\
    \frac{\xi+2\sqrt{\xi}+2}{2+\sqrt{\xi}}\beta_p-\alpha_p^{-1}\ln\frac{\sqrt{1+\sqrt{\xi}}}{2+\sqrt{\xi}},
                & \hbox{ if }~\xi>2.
  \end{array}
\right.
\]
Let $\Lambda(z)=\exp(-e^{-z}), ~ z\in \mathbb{R},$ be the distribution function of the Gumbel distribution or type I extreme-value distribution. Its probability density function is $\lambda(z)=\exp(-z-e^{-z})$.

\begin{theorem}\label{main1} We have
\begin{equation}\label{Gumbel}
A_{p,\xi}\Big(\max_{1\le i<j\le
p}l_{ij}-B_{p,\xi}\Big)
\end{equation}
converges weakly to a probability distribution with distribution  function $G_{\xi}(z)$,
where
$G_\xi(z)=\Lambda(z)$ for $\xi\in [0, 2]$,
and
\begin{equation}\label{G-xi}
G_{\xi}(z)=\int^{\infty}_{-\infty}\exp\Big(-e^{(y-z)/\eta}- \eta
e^{-z}\int^{e^{(y-z)/\eta}}_0s^{-1-\eta}(1-e^{-s})ds
\Big)\lambda(y)dy, ~~ z\in \mathbb{R},
\end{equation}
for $\xi>2$, where
$\eta=(2+\sqrt{\xi})(\xi+\sqrt{\xi})^{-1}$.
\end{theorem}

Evidently,  $G_{\xi}(z)$ from \eqref{G-xi} is a continuous function and is increasing in $z\in \mathbb{R}$. Also, $G_{\xi}(-\infty)=0$ and $G_{\xi}(\infty)=1.$ Hence $G_{\xi}(\cdot)$ is a probability distribution function. It is tantalizing to have an interpretation of  $G_{\xi}(z)$. At the first sight the function $G_{\xi}(z)$ looks like a convolution of function $\lambda(y)$ and another one. However, we stop trying the thought because of the term $e^{-z}$ in the integral from \eqref{G-xi} while other terms are functions of $z-y$.

Obviously, the class of distribution functions $G_{\xi}(\cdot)$  has a transition at $\xi=2$. A quick glimpse over  the constant $\ln[\sqrt{2}\pi/\arcsin(\sqrt{2}-1)]$ in the definition of $B_{p,\xi}$ seems a good indication.  As $\xi \leq 2$, the limiting distribution is the classical Gumbel distribution, however, $G_{\xi}(z)$ is  a new one as $\xi>2$ to the best of our knowledge.

It is straightforward to check that Theorem \ref{main1} for the case $\xi=2$ is the same as \eqref{cozy} at $m=2$.

Theorem \ref{main1} does not hold when the Gaussian assumption is removed. For example, assume $\bz_{ii}=0$ for each $i$ in the  matrix $\bd{Z}=(\bd{z}_{ij})$ and $\bz_{ij}$ follows the uniform distribution on $[-2, 2]$ for any $i<j$. This corresponds to $\xi=0$ in Theorem \ref{main1}. Then $l_{ij}=|\bd{z}_{ij}|$ and therefore $\{l_{ij};\, 1\leq i<j \leq p\}$ are i.i.d. random variables with the uniform distribution on $[0, 2]$. In this case, one can easily check  that $(p^2/4)(2-\max_{1\le i<j\le
p}l_{ij})$ converges weakly to the unit exponential distribution Exp($1$).

In principle, one can generalize Theorem \ref{main1} to $3\times 3$ and $4\times 4$ minors. An encouraging sign is that, similar to \eqref{Dimsum}, the eigenvalues of $3\times 3$ and $4\times 4$ matrices are the roots of their cubic and quartic characteristic polynomials, respectively. As a result these eigenvalues have explicit formulas, but this is not true for $k\times k$ minors with $k\geq 5$ by the well-known works of Abel and Galois. A drawback is that the Stein Poisson approximation method, as the initial step to derive the limiting distribution in Theorem \ref{main1}, may not be easily adjusted to modify the $k\times k$ minors for $k=3,4$. The  $2\times 2$ minors have simpler structures after all. Thus, in order to study the $k\times k$ minors with $k\geq 5$, one has to discover a new method.

Look at the distributions $G_{\xi}(\cdot)$ for $0\leq \xi \leq 2$ and $\xi>2$, they appear very differently. One may wonder if the distribution from \eqref{G-xi} is special to the $2\times 2$ minors. In other words, for the case $\xi > 2$, do we always expect the distribution from \eqref{G-xi} for arbitrary $k\geq 2$? Or there exists a threshold $k_0$ such that the limiting distribution is $\Lambda(z)=\exp(-e^{-z})$ for any $k\times k$ minors with $k\geq k_0$, regardless of the value of $\xi \in [0, \infty)$?

Finally let us explain the method we employed to derive the Gumbel distribution and that from \eqref{G-xi}. In fact, reviewing $l_{ij}$ from \eqref{Dimsum}, conditioning on the diagonal entries $\bz_{ii}$'s, we apply the Stein Poisson approximation method. Via this step we reduce the study of the distribution function of  $\max_{1\le i<j\le p}l_{ij}$ to the moment generating function of a $U$-statistic $\sum_{1\leq i<j \leq p}\psi(\bar\bz_i, \bar\bz_j)$, where $\bar\bz_i$'s are i.i.d. non-Gaussian random variables with finite upper bounds. Second, by the H\'ajek projection, the $U$-statistic is reduced to the sum of some i.i.d. random variables up to certain errors that are under control. Then we study the moment generating functions of the new random variables.

Now we would like to comment on the deformed GOE in a different direction: the properties of the eigenvalues of the deformed GOE $\bd{Z}=\bd{Z}_{\xi}$ itself instead of those of minors of $\bZ$.  In fact, replacing the diagonal entries with $0$ will not affect the long-term behavior of its largest eigenvalue and empirical spectral distribution of $\bd{Z}$. Precisely,  the largest eigenvalue of $\bd{Z}$ converges to the Tracy-Widom law (see \cite{kho12}, \cite{Ru2006} or \cite{Shasha}) as $p\to \infty$. The empirical spectral distribution of $\bd{Z}$ converges to the semi-circle law by Lidskii’s theorem (see \cite{Lid1950} or Lemma 2.3 from \cite{Bai1999}).

Obviously  the maximum of the largest eigenvalues of all $1\times 1$ principal minors is simply the maximum of $p$ i.i.d. $N(0, \xi)$-distributed random variables. Consequently, the maximum has an asymptotic Gumbel distribution regardless of the values of $\xi$ as long as $\xi \geq 0$. According to the discussion in the above paragraph, we can regard \eqref{G-xi} as an interpolation between the Gumbel distribution and the Tracy-Widom distribution since the largest eigenvalues of $\bd{Z}$ and its $p\times p$ principal minor are the same.


It is time to study our original problem on the Wishart case. Recall \eqref{NEWS}. Now we consider $T_{n,p,m}$ with $m=2$. We assume $p$ changes with $n$ for mathematical rigor.

%
%
%
%
%
%
%
%
%
%
%

\begin{theorem}\label{main0}  Assume $E(x_{11}^6)<\infty$. Set $\xi=\mbox{Var}(x_{11}^2) \in [0, \infty).$ If $p=p_n\to\infty$ and
$p=o(n^{1/7})$ as $n\to\infty$, then
\[
A_{p,\xi}\left(\frac{T_{n,p,2}-n}{\sqrt{n}}-B_{p,\xi}\right)
\]
converges weakly to a probability distribution with distribution  function $G_{\xi}(z)$, where $A_{p,\xi}$, $B_{p,\xi}$ and
$G_{\xi}(z)$ are defined in Theorem \ref{main1}.
\end{theorem}

A special case of the above theorem is that $\{x_{ij}\}$ are i.i.d. symmetric Bernoulli random variables, that is, $P(\xi_{11}=\pm 1)=1/2.$ In this case, $\xi_{11}$ has mean zero, variance one and $\xi=\mbox{Var}(x_{11}^2)=0$.
Since $\mathbf{X}^{T}\mathbf{X}=(w_{ij})_{1\le i,j\le
p}$, we know $w_{ii}=n$ and $T_{n,p,2}-n=\max_{1\leq i< j\leq p}|w_{ij}|$. Our result here is consistent with Theorem 3 and Remark 2.1 from \cite{CJ2011} for the statistic $J_n$ there, which handles the maximum of the absolute values of the non-diagonal entries of a sample covariance matrix. Relevant work in the literature include \cite{CFJ131988, CJ12, FSZ18, Jiang04, JP24, LLR10, LQR12, LR06, SZ14, Zhou07}, which is by no means exhaustive.

The proof of Theorem \ref{main0} is based on Theorem \ref{main1}. In fact, a high-dimensional central limit theorem (see Theorem \ref{normalapp}) bridges the sample covariance matrix $\bd{W}$ and the deformed GOE $\bd{Z}$ discussed in \eqref{chifan}. The property of $\bd{Z}$ appeared in Theorem \ref{main1} then carries over to that of $\bd{W}$ in Theorem \ref{main0}. This not only confirms the conjecture stated between \eqref{callback} and \eqref{cozy} that $(T_{n,p, m}-n)/\sqrt{n}$ and $\tilde{T}_{p, m}$ have similar behavior if $\xi=\mbox{Var}(x_{11}^2) \leq 2$ but also discovers that the similar behavior holds for $\xi> 2.$

As mentioned above, we rely on the limiting distribution of $\tilde{T}_{p, 2}$ in \eqref{NEWS1} for the deformed GOE to get that of $T_{n,p,2}$. By the same argument, in order to get the asymptotic distribution of $T_{n,p,m}$ for $m\geq 3$, one needs to work out $\tilde{T}_{p, m}$ first. Recall the discussion on
 $\tilde{T}_{p, m}$ for $m\geq 3$ below Theorem \ref{main1}. We leave this as a future problem.


\section{The proofs of main results}\label{proof:main}

In this section, we will prove Theorems~\ref{main1} and \ref{main0}. The proofs are quite involved, and therefore for clarity we will provide a  few of lemmas first without proofs. Their proofs will be given separately in Section~\ref{lemmas}.

\subsection{Proof of Theorem~\ref{main1}}

We will continue to use the notation $\xi, \bz_{ij}, \bZ_{\{i,j\}}, l_{ij}$
 and $\bZ=(\bz_{ij})_{1\le i,j\le p}$  given at the beginning of Section \ref{mainresults}.
Let $\phi$ and $\Phi$ denote the density function
and cumulative distribution function of $N(0, 1)$ respectively, that is,
\[
\phi(x)=\frac{1}{ \sqrt{2\pi}}e^{-x^2/2} ~~\mbox{ and }~~
\Phi(x)=\int^x_{-\infty}\phi(s)ds,~~~ x \in \mathbb{R}.
\]


\noindent For any real number $c$,  define a density function $\phi_c$ with
\begin{equation}\label{tongzhia}
\phi_{c}(x):=\left\{
            \begin{array}{ll}
              \frac{\phi(x)}{\Phi(c)}, & x\le c; \\
              0, & x>c.
            \end{array}
          \right.
\end{equation}
If $\bz \sim N(0,1)$,  clearly,  $\phi_{c}$ is the conditional density function of $\bz$ given $\bz\le c$.
Recall $\alpha_p=\sqrt{2\ln p}$ as defined in the beginning of Section~\ref{mainresults}. Set
\begin{equation}\label{standard}
\beta_p=\alpha_p-\alpha_p^{-1}\ln\big(\sqrt{2\pi}\alpha_p\big)~~\mbox{
and}~~~c_p(y)=\beta_p+\alpha_p^{-1}y, ~~~y\in \mathbb{R}.
\end{equation}
For brevity, we now explain where the constants $\alpha_p, \beta_p$ and $c_p(y)$ come from and this also serves the purpose of analysis in the future. In fact, in the classic extreme-value theory,  it is well known that
$\alpha_p(\max\limits_{1\le j\le p}\bz_j-\beta_p)$ converges weakly to the Gumbel distribution $\Lambda(y)=\exp(-e^y)$, $y\in \mathbb{R}$, or equivalently,
\begin{equation}\label{EV}
P\Big(\alpha_p\big(\max\limits_{1\le j\le p}\bz_j-\beta_p\big)\le
y\Big)=\big[\Phi(c_p(y))\big]^p=\Big[1-\frac{[1+o(1)]e^{-y}}{p}\Big]^p\to
\Lambda(y)
\end{equation}
for any $y\in \mathbb{R}$.  This can be easily verified by
using the following well-known fact that
\begin{equation}\label{GaussianTail}
\Phi(-z)=1-\Phi(z)=[1+o(1)]\frac{\phi(z)}{z}=[1+o(1)]\frac{1}{\sqrt{2\pi}z}\exp\Big(-\frac12z^2\Big)
\end{equation}
as $z\to\infty$.

First, notice the case $\xi=0$ is trivial. In fact, in virtue to \eqref{Dimsum}, $l_{ij}=|\bz_{ij}|$ for all $i\ne j$. Hence $\max_{1\leq i<j\leq p}l_{ij}=\max_{1\leq i<j\leq p}|\bz_{ij}|$, a maximum of $p(p-1)/2$ i.i.d. random variables. Now $\alpha_p=\sqrt{2\ln p}$, $A_{p,0}=\sqrt{2}\alpha_p$ and $B_{p,0}=\sqrt{2}\alpha_p-(\sqrt{2}\alpha_p)^{-1}\ln (\sqrt{2}\alpha_p\sqrt{2\pi}).$ Set $A_p=A_{p,0}$ and $B_p=B_{p,0}.$ Then, by independence and an argument similar to \eqref{EV} and \eqref{GaussianTail}, we obtain
\begin{eqnarray*}
P\Big(A_{p,0}\big(\max_{1\le i<j\le
p}l_{ij}-B_{p,0}\big)\le z\Big)
&=&  P\big(|\bz_{12}|\leq B_p+A_p^{-1}z\big)^{p(p-1)/2} \\
&=& \big[1-2\Phi\big(-B_p-A_p^{-1}z\big)\big]^{p(p-1)/2}\\
& \to & \exp(-e^{-z})
\end{eqnarray*}
for any $z\in \mathbb{R}.$

So we always assume $\xi>0$ in the remaining proof.

Now, for convenience, we write $\bz_i=\bz_{ii}/\sqrt{\xi}$. Hence, $\bz_i$'s are i.i.d. $N(0, 1)$-distributed random variables and are independent of $\bz_{ij}$'s.
With this change of notation, we now have
\begin{equation}\label{mooncake}
l_{ij}=\lambda_1(\bZ_{\{i,j\}})=
\sqrt{\bz_{ij}^2+\frac{\xi}{4}(\bz_i-\bz_j)^2}+
\frac{\sqrt{\xi}}{2}(\bz_i+\bz_j)
\end{equation}
 and the $2\times 2$ principal minor
\[
\bZ_{\{i,j\}}=
\left(
          \begin{array}{cc}
            \sqrt{\xi}\bz_i & \bz_{ij} \\
            \bz_{ij} & \sqrt{\xi}\bz_j \\
          \end{array}
        \right)
\]
for all $1\le i< j\le p.$ Assume $\{t_p\}$ is a sequence of positive numbers such that
\begin{equation}\label{tp}
\lim\limits_{p\to\infty}\frac{t_p}{\sqrt{2\ln p}}:=r\in \left(\sqrt{\xi},
\sqrt{\xi}+\frac{2}{\sqrt{\xi}}\right).
\end{equation}
Our objective is to estimate the probability $P(\max_{1\le i<j\le
p}l_{ij}\le t_p)$.   To this end, for all real numbers $x$,
$y$ and $t$, define
\begin{equation}\label{q(x,y)}
q(x,y;t)=P\Big(\bz_{12}^2>(t-\sqrt{\xi}x)(t-\sqrt{\xi}y)\Big),
\end{equation}
which is equal to
$2[1-\Phi(\sqrt{(t-\sqrt{\xi}x)(t-\sqrt{\xi}y)})]$ if
$(t-\sqrt{\xi}x)(t-\sqrt{\xi}y)>0$, and $1$ otherwise.

\vspace{10pt}

In the following we will study
$P(\max_{1\le i<j\le p}l_{ij}\le t_p)$. The analysis is divided into Lemmas \ref{lem0}-\ref{estofq}. Their proofs are presented in Section \ref{lemmas}.

\begin{lemma}\label{lem0} Assume condition
\eqref{tp} holds. Then
\begin{equation*}
\lim\limits_{p\to\infty}\left|P\left(\max_{1\le i<j\le p}l_{ij}\le
t_p\right)-E\exp\left(-\sum_{1\le i<j\le p}q(\bz_i,\bz_j;t_p)\right)\right|=0.
\end{equation*}
\end{lemma}

\begin{lemma}\label{rep0}
Let $f_p(z_1, \cdots, z_p)$ be a sequence of symmetric
functions of $z_1,\cdots, z_p$, and for some constant $b>0$, $|f_p(z_1,
\cdots, z_p)|\le b$ for all $(z_1, \cdots, z_p)\in \mathbb{R}^p$ and $p\ge 2$.  Then
\[
Ef_p(\bz_1, \cdots, \bz_{p-1},
\bz_p)-\int^\infty_{-\infty}Ef_p(\bar\bz_1, \cdots,
\bar\bz_{p-1}, c_p(y))\lambda(y)dy\to 0
\]
as $p\to\infty$, where for each $y\in \mathbb{R}$, $\bar\bz_1, \cdots,
\bar\bz_{p-1}$ are i.i.d. random variables with density function
$\phi_{c_p(y)}(\cdot)$  and $\lambda(y)=\exp(-y-e^{-y})$.
\end{lemma}

By combining Lemmas~\ref{lem0} and \ref{rep0} we immediately have
the following lemma.

\begin{lemma}\label{rep1} Under condition \eqref{tp},
\[
\lim\limits_{p\to\infty}\left|P\left(\max_{1\le i<j\le p}l_{ij}\le t_p\right)
-\int^\infty_{-\infty}Ee^{-Q_p-U_p}\lambda(y)dy\right|=0,
\]
where $Q_p:=\sum_{1\le i<j\le p-1}q(\bar\bz_i, \bar\bz_j;t_p)$,
$U_p:=\sum^{p-1}_{j=1}q(c_p(y), \bar\bz_j;t_p)$, and $\bar\bz_1,
\cdots, \bar\bz_{p-1}$, as defined in Lemma~\ref{rep0}, are i.i.d.
random variables with density function $\phi_{c_p(y)}(\cdot)$ for each
$y$.  Furthermore, if $Ee^{-Q_p-U_p}$ converges to $g(y)$ for almost every
$y\in \mathbb{R}$ as $p\to\infty$, then it follows from the dominated convergence theorem that
\[
\lim\limits_{p\to\infty}P\left(\max_{1\le i<j\le p}l_{ij}\le
t_p\right)=\int^\infty_{-\infty}g(y)\lambda(y)dy.
\]
\end{lemma}


From  Lemma~\ref{rep1},  we need to handle some functions of independent random variables $\bar\bz_1,
\cdots, \bar\bz_{p-1}$.  Note that these random variables have a
density function $\phi_{c_p(y)}(\cdot)$ for each $y\in \mathbb{R}$,
where $c_p(y)$ is defined in \eqref{standard} with property
\[
\lim\limits_{p\to\infty}\frac{c_p(y)}{\sqrt{2\ln p}}=1,~~ y\in
\mathbb{R}.
\]
Since we only need to find the limit of the
expectation for each $y$ as stated in Lemma~\ref{rep1}, the symbol ``$y$" is suppressed from now on for
simplicity.  In fact, sometimes we only need $c_p$ to be a sequence of
 real numbers satisfying
\begin{equation}\label{cp}
\lim\limits_{p\to\infty}\frac{c_p}{\sqrt{2\ln p}}=1.
\end{equation}
There is no need to assume $c_p$ has an explicit or close form as in \eqref{standard}. Condition \eqref{cp} will be specified when this is the case.

As in Lemma \ref{rep1}, $\bar\bz_1, \cdots, \bar\bz_{p-1}$ are independent
random variables with density function $\phi_{c_p(y)}(\cdot)$. We will investigate
the limit of $Ee^{-Q_p-U_p}$, where
\begin{equation}\label{QU}
Q_p=\sum_{1\le i<j\le p-1}q(\bar\bz_i, \bar\bz_j;t_p),~~
U_p=\sum^{p-1}_{j=1}q(c_p, \bar\bz_j;t_p).
\end{equation}

Review that $q(x,y;t)$ is defined in \eqref{q(x,y)}.  Under
conditions \eqref{tp} and \eqref{cp}, we have, for all large $p$,
\begin{equation}\label{q(x,y;tp)}
q(x,y;t_p)=2\left[1-\Phi\left(\sqrt{\left(t_p-\sqrt{\xi}x\right)\left(t_p-\sqrt{\xi}y\right)}\,\right)\right]
\end{equation}
for $x, y\le c_p$ since $t_p-\sqrt{\xi}c_p>0$ for all large $p$.
Also note that
$\bar\bz_j\le c_p$ for all $j$.   Hence, we have
\[
q(\bar\bz_i,\bar\bz_j;
t_p)=2\left[1-\Phi\left(\sqrt{\left(t_p-\sqrt{\xi}\bar\bz_i\right)\left(t_p-\sqrt{\xi}\bar\bz_j\right)}\,\right)\right].
\]
Now we define
\begin{equation}\label{q(x)}
q(x;t_p)=Eq(x,\bar\bz_2;t_p),~~~x\le c_p,
\end{equation}
and
\begin{equation}\label{qtp}
q(t_p)=Eq(\bar\bz_1,\bar\bz_2; t_p)=
Eq(\bar\bz_1;t_p).
\end{equation}



In the following lemma, for clarity we need to mention that $p^{\alpha}q(c_p;t_p)$ stands for $p^{\alpha}\cdot q(c_p;t_p)$ for any $\alpha\in \mathbb{R}$.

\begin{lemma}\label{rep2} Assume conditions \eqref{tp} and \eqref{cp} hold.\\
(a). If $pq(c_p;t_p)$ and $p^2q(t_p)$ are bounded uniformly in
$p$, then, as $p\to\infty$,
\begin{equation}\label{rep2a}
Ee^{-Q_p-U_p}-E\exp\Big(d_p-(p-2)\sum^{p-1}_{j=1}q(\bar\bz_j;t_p)\Big)
 \to 0,
\end{equation}
where $Q_p$ and $U_p$ are given in \eqref{QU}, and
$d_p=\frac12(p-1)(p-2)q(t_p)-(p-1)q(c_p;t_p)$. In particular, we have the following two conclusions.\\
(b). If $\lim\limits_{p\to\infty}pq(c_p;t_p)=0$ and
$\lim\limits_{p\to\infty}(1/2)p^2q(t_p)=\gamma\in (0,\infty)$, then
\begin{equation}\label{rep2b}
\lim\limits_{p\to\infty}Ee^{-Q_p-U_p}=e^{-\gamma}.
\end{equation}
(c). Assume $\lim\limits_{p\to\infty}pq(c_p;t_p)=\tau\in (0,\infty)$ and
$\lim\limits_{p\to\infty}(1/2)p^2q(t_p)=\gamma_1\in (0,\infty)$. If
\[
\lim\limits_{p\to\infty}p^{j+1}E\big(q(\bar\bz_1;t_p)^j\big):=\gamma_j\in
[0,\infty)
\]
for each $j\ge 2$, then
\begin{equation}\label{irregularlimit}
\lim\limits_{p\to\infty}Ee^{-Q_p-U_p}=\exp\left(-\tau-\sum^\infty_{j=1}(-1)^{j-1}\frac{\gamma_j}{j!}\right).
\end{equation}
\end{lemma}
The idea of deriving \eqref{rep2a} is actually inspired by the H\'ajek projection of $U$-statistics. The projection  reduces a $U$-statistic to the sum of independent random variables and small terms. Through this one usually  gets the central limit theorem of the $U$ statistic; see, for example, \cite{Vaart98}. Fortunately this philosophy  works in our setting for the moment generating functions of $U$-statistics, too.

In Lemmas~\ref{lem2} and \ref{estofq} next, we will estimate $q(x,t_p)$  and $q(t_p)$, respectively.

\begin{lemma}\label{lem2}  Assume conditions \eqref{tp} and \eqref{cp} hold. Then the following are true for any $\xi>0$.\\
(i). We have
\begin{equation}\label{q(x)-upper}
q(x;t_p)\le \frac{2[1+o(1)]}{\sqrt{2\pi}\sqrt{t_p-\sqrt{\xi}x}}
\frac{\Phi(\bar
c_p)}{\sqrt{t_p-\sqrt{\xi}c_p}}\exp\left(-\frac12\left(\frac{4-\xi}{4}t_p^2-\frac{\xi^2}{4}x^2-\frac{(2-\xi)\sqrt{\xi}}{2}xt_p\right)
\right)
\end{equation}
uniformly over $x\le c_p$ as $p\to\infty$, where $\bar c_p:=
\bar c_p(t_p,x)=c_p-\frac12\sqrt{\xi}(t_p-\sqrt{\xi}x)$.
\\
(ii). For any sequence of positive numbers $b_p$ with
$b_p\to\infty$ and $b_p/c_p\to 0$ as $p\to\infty$, we have
\begin{equation}\label{q(x)-approx}
q(x;t_p)=\frac{2[1+o(1)]\exp\left(-\frac12\Big(\frac{4-\xi}{4}t_p^2-\frac{\xi^2}{4}x^2-\frac{(2-\xi)\sqrt{\xi}}{2}xt_p\Big)
\right)}{\sqrt{2\pi}\sqrt{\big(t_p-\sqrt{\xi}x\big)\left[(1-0.5\xi)t_p+0.5\xi^{3/2}x\right]}}
\end{equation}
 uniformly over $x_p\le x\le c_p$ as $p\to\infty$, where
$x_p=\frac1{\sqrt{\xi}}\big(t_p-\frac{2c_p}{\sqrt{\xi}}\big)+b_p$.\\
(iii). Assume $b_p$ is as in (ii). Then
\begin{equation}\label{q(x)-approx-small-x}
q(x;t_p)=
\frac{2[1+o(1)]\exp(-\frac12c_p^2)}{\pi\sqrt{t_p-\sqrt{\xi}c_p}}\cdot
\frac{\exp\big(-\frac{1}{2}(t_p-\sqrt{\xi}x)(t_p-\sqrt{\xi}c_p)\big)}{\sqrt{t_p-\sqrt{\xi}x}\cdot\big[\sqrt{\xi}(t_p-\sqrt{\xi}x)-2c_p\big]}
\end{equation}
uniformly over $x\le \bar x_p=
\frac1{\sqrt{\xi}}\big(t_p-\frac{2c_p}{\sqrt{\xi}}\big)-b_p$ as
$p\to\infty$.

\end{lemma}




\begin{lemma}\label{estofq}  Assume \eqref{tp} and \eqref{cp}. Recall \eqref{qtp}. The following hold as $p\to\infty$.\\
\noindent (a). If $\xi\in (0,2)$, then
\begin{equation}\label{case1}
q(t_p)=[1+o(1)]\sqrt{\frac{2(2+\xi)}{\pi(2-\xi)}}\cdot\frac{\exp\Big(-\frac{t_p^2}{2+\xi}\Big)}{t_p}.
\end{equation}


\noindent (b). If $\xi=2$, recalling $r$ in \eqref{tp}, then
\begin{equation}\label{case2}
 q(t_p)=\frac{[1+o(1)]\exp\big(-t_p^2/4\big)}{\sqrt{2}\pi}\int^{\sqrt{2}/r}_{1-\sqrt{2}/r}s^{-1/2}(1-s)^{-1/2}ds.
 \end{equation}

\noindent (c). If $\xi\in (2,\infty)$,  then
\begin{equation}\label{case3}
q(t_p)=\frac{8[1+o(1)]q(c_p;t_p)\phi(c_p)}{(\xi^2-4)c_p-(\xi-2)\sqrt{\xi}t_p}.
\end{equation}
Furthermore,  for $j\ge 2$, we have
\begin{equation}\label{j-moments}
E\left(q(\bar\bz_{1};t_p)^j\right)=\frac{4[1+o(1)]q(c_p;t_p)^{j}\phi(c_p)}{(j\xi^2-4)c_p-j(\xi-2)\sqrt{\xi}t_p}
+o\Big(q(c_p;t_p)^{j-1}q(t_p)\Big).
\end{equation}
\end{lemma}

Some basic but important calculations are needed in the discussions ahead. They are collected in the following two lemmas.

\begin{lemma}\label{chores1} Let $\alpha_p, \beta_p, c_p$ be as in \eqref{standard}. Let  $q(x; t_p)$, $q(t_p)$ and $q(\bar\bz_{1};t_p)$ be as in Lemmas \ref{lem2} and \ref{estofq}, respectively.  Recall $A_{p,\xi}$ and $B_{p,\xi}$  in Theorem \ref{main1}. For fixed  $z\in \mathbb{R}$, set
$t_p=B_{p,\xi}+A_{p,\xi}^{-1}z$. The following hold. \\
(i) $\phi(c_p)/c_p=\big(1+o(1)\big) e^{-y}/p$, where $\phi(\cdot)$ is the density function of $N(0, 1).$\\
(ii) If $\xi \in (0, 2]$, then
$\lim\limits_{p\to\infty}p^2q(t_p)=2e^{-z}.$\\
(iii) If $\xi \in (0, 2]$, then $q(c_p;t_p)\phi(c_p)= o(q(t_p)).$
\end{lemma}

\begin{lemma}\label{chores2} Let $\alpha_p, \beta_p, c_p$ be as in \eqref{standard}. Let  $q(x; t_p)$, $q(t_p)$,  $q(\bar\bz_{1};t_p)$, and $t_p$ be as in Lemmas \ref{lem2},  \ref{estofq} and \ref{chores1}, respectively. Assume $\xi>2.$ Set $\eta=(2+\sqrt{\xi})(\xi+\sqrt{\xi})^{-1}$. Then the  following are true.\\
(a) $(j\xi^2-4)c_p-j(\xi-2)\sqrt{\xi}t_p=
4\big(1+o(1)\big)\eta^{-1}(j-\eta)c_p$ for each integer $j\geq 0$. \\
(b) $
\lim\limits_{p\to\infty}pq(c_p;t_p)=\exp\big(\eta^{-1}(y-z)\big).$\\
(c) $\lim_{p\to\infty}p^2q(t_p)=2\eta\tau e^{-y}/(1-\eta),$ where $\tau=e^{(y-z)/\eta}$.\\
(d) $\lim_{p\to\infty}p^{j+1}E\big(q(\bar\bz_{1};t_p)^j\big)=
\eta(j-\eta)^{-1} \tau^je^{-y}$ for $j\ge 2$.
\end{lemma}

\medskip

With the preparation above we now proceed to prove Theorem~\ref{main1}.
To show \eqref{Gumbel}, we need to verify
\begin{equation*}
\lim\limits_{p\to\infty}P\Big(A_{p,\xi}\big(\max_{1\le i<j\le
p}l_{ij}-B_{p,\xi}\big)\le z\Big)=G_{\xi}(z)
\end{equation*}
for any $z\in \mathbb{R}$. Notice
\begin{eqnarray}\label{diyige}
 P\Big(A_{p,\xi}\big(\max_{1\le i<j\le
p}l_{ij}-B_{p,\xi}\big)\le z\Big)=P\left(\max_{1\le i<j\le
p}l_{ij} \le t_p\right)
\end{eqnarray}
where $t_p:=B_{p,\xi}+A_{p,\xi}^{-1}z$. Then, by the definition of $r$ in \eqref{tp}, we know
\begin{eqnarray*}
r=\lim_{p\to\infty}\frac{t_p}{\sqrt{2\ln  p}}=\lim\limits_{p\to\infty}\frac{t_p}{c_p} =\left\{
 \begin{array}{ll}
 \sqrt{2+\xi}, & \xi\in (0,2]; \\
 \frac{\xi+2\sqrt{\xi}+2}{2+\sqrt{\xi}}, & \xi> 2.
 \end{array}
 \right.
\end{eqnarray*}
It is trivial to see $r\in (\sqrt{\xi},
\sqrt{\xi}+(2/\sqrt{\xi}))$ for every $\xi>0$. Hence the requirement in \eqref{tp} is satisfied.  For every $y\in \mathbb{R}$, let $c_p=c_p(y)$ be as in \eqref{standard}. We know $c_p/\alpha_p=c_p/\sqrt{2\ln p} \to 1$ as $p\to\infty$, which confirms \eqref{cp}. All requirements appeared in Lemmas \ref{lem0}-\ref{rep2} are satisfied due to the assumptions \eqref{tp} and \eqref{cp}. Therefore they
all hold true.  Next, to prove the theorem, we consider the cases $\xi \in (0, 2]$ and $\xi>2$ separately.

{\it The case $\xi \in (0, 2]$}. Observe
\begin{eqnarray*}
  pq(c_p;t_p)
  =\frac{c_p}{p\phi(c_p)}\cdot p^2q(t_p)\cdot\frac{q(c_p;t_p)\phi(c_p)}{q(t_p)}\cdot \frac{1}{c_p}=o\left(\frac{1}{c_p}\right)=o(1)
\end{eqnarray*}
by the three assertions from Lemma \ref{chores1} and the fact $c_p\sim \sqrt{2\ln p}$. Also, by Lemma \ref{chores1}(ii), $\frac12p^2q(t_p)=\gamma\in (0,\infty)$ with $\gamma=e^{-z}$. We hence have from Lemma \ref{rep2}(b) that
\begin{eqnarray*}
    \lim_{p\to\infty}Ee^{-Q_p-U_p}=e^{-e^{-z}}.
\end{eqnarray*}
Recall $z\in \mathbb{R}$ is fixed in advance. Take $g(y) \equiv e^{-e^{-z}}$ for each $y \in \mathbb{R}.$ We obtain from Lemma \ref{rep1} that
\[
\lim_{p\to\infty}P\left(\max_{1\le i<j\le p}l_{ij}\le
t_p\right)=\int^\infty_{-\infty}g(y)\lambda(y)dy=e^{-e^{-z}}
\]
since $\lambda(y)$ is a probability density function. This and \eqref{diyige} conclude \eqref{Gumbel} for the case  $\xi \in (0, 2]$.

{\it The case $\xi >2$}. Recall $\eta=(2+\sqrt{\xi})(\xi+\sqrt{\xi})^{-1}$. From Lemma \ref{chores2} we see
\begin{eqnarray*}
&&
\lim\limits_{p\to\infty}pq(c_p;t_p)=\tau:=\exp\big(\eta^{-1}(y-z)\big);\\
&&\lim_{p\to\infty}\frac{1}{2}p^2q(t_p)=\gamma_1=\frac{\eta\tau e^{-y}}{1-\eta};\\
&& \lim_{p\to\infty}p^{j+1}E\big(q(\bar\bz_{1};t_p)^j\big)=\gamma_j:=
\eta(j-\eta)^{-1} \tau^je^{-y}, ~~  j\ge 2.
\end{eqnarray*}
By Lemma \ref{rep2}(c),
\begin{equation*}
\lim_{p\to\infty}Ee^{-Q_p-U_p}=\exp\left(-\tau-\sum^\infty_{j=1}(-1)^{j-1}\frac{\gamma_j}{j!}\right).
\end{equation*}
Now we evaluate the above infinite series. In fact,
\begin{eqnarray*}
\sum^\infty_{j=1}(-1)^{j-1}\frac{\gamma_j}{j!}&=&\eta
e^{-y}\tau^{\eta}\sum^\infty_{j=1}\frac{(-1)^{j-1}}{j!}\frac{\tau^{j-\eta}}{j-\eta}\\
&=&\eta
e^{-y}\tau^{\eta}\sum^\infty_{j=1}\frac{(-1)^{j-1}}{j!}\int^{\tau}_0s^{j-1-\eta}ds\\
&=&-\eta
e^{-y}\tau^{\eta}\int^{\tau}_0s^{-1-\eta}\sum^\infty_{j=1}\frac{(-s)^{j}}{j!}ds\\
&=&\eta
e^{-y}\tau^{\eta}\int^{\tau}_0s^{-1-\eta}(1-e^{-s})ds\\
&=&\eta e^{-z}\int^{\tau}_0s^{-1-\eta}(1-e^{-s})ds,
\end{eqnarray*}
where in the last step, the formula $\tau:=\exp((y-z)/\eta)$ is used. Consequently,
\begin{eqnarray*}
\lim_{p\to\infty}Ee^{-Q_p-U_p}=\exp\left(-\exp\big(\eta^{-1}(y-z)\big)-\eta e^{-z}\int^{e^{(y-z)/\tau}}_0s^{-1-\eta}(1-e^{-s})ds\right).
\end{eqnarray*}
We conclude from Lemma \ref{rep1} that
\begin{eqnarray*}
&&\lim_{p\to\infty}P\left(\max_{1\le i<j\le p}l_{ij}\le
t_p\right)\\
&=& \int^\infty_{-\infty}\exp\left(-e^{(y-z)/\eta}-\eta e^{-z}\int^{e^{(y-z)/\eta}}_0s^{-1-\eta}(1-e^{-s})ds\right)\lambda(y)dy
\end{eqnarray*}
where $\lambda(y)=\exp(-y-e^{-y})$. This together with \eqref{diyige} yields \eqref{G-xi}. The proof is completed. \eop

\subsection{Proof of Theorem~\ref{main0}}

We first cite a high-dimensional central limit theorem.  The following lemma is a special case of Theorem 1.1 from
Rai\u{c}~\cite{Raic19}.

\begin{lemma}\label{normalapp}
Let $\bY_k$, $k=1,\ldots, n,$ be $n$ independent $\mathbb{R}^d$-valued random
vectors with mean $E(\bY_k)=0$ and covariance matrix
$\mbox{Cov}(\bY_k)=I_d$, where $I_d$ is the $d\times d$ identity matrix.
Assume $\bT$ is a $d$-dimensional Gaussian random vector with zero
mean and $\mbox{Cov}(\bT)=I_d$. Then
\[
\sup_{A\in \cal C}\left|P\left(\frac{1}{\sqrt{n}}\sum^n_{k=1}\bY_k\in
A\right)-P\left(\bT\in A\right)\right|\le
(42d^{1/4}+16)n^{-1/2}E\left(|\bY_1|_2^{3}\right),
\]
where $\cal C$ denotes the collection of all measurable convex sets
$A\subset \mathbb{R}^d$, and $|\cdot|_2$ denotes the Euclidean norm for
vectors in $\mathbb{R}^d$.
\end{lemma}

Recall $\mathbf{X}=(x_{ij})_{1\le i\le
n,1\le j\le p}$, where the $np$ random variables $\{x_{ij}\}$ are i.i.d. and each has  mean $0$ and variance $1.$
The (white) Wishart matrix is defined by $\mathbf{W}=(w_{ij})_{p\times p}=\mathbf{X}^{T}\mathbf{X}$. Obviously, $w_{ij}=\sum_{i=1}^nx_{ki}x_{kj}.$
For each $1\le k\le n$, write
$y_{k:ii}=x_{ki}^2-1$ and
$y_{k:ij}=x_{ki}x_{kj}$ for $1\le i\ne j\le p$, and define a vector
\[
\bY_k=\big(y_{k:11},\ldots, y_{k:pp}, y_{k:ij}; 1\le i<j\le p\big)^{T},
\]
that is,  for each $k$, $\bY_k$ is a random vector whose coordinates consist of all upper triangular entries of the symmetric matrix $\big(y_{k:ij}\big)_{1\le i,j\le p}$ as well as the $p$ diagonal entries of $\mathbf{W} - \mathbf{I}_p.$ Let us prove the theorem by differentiating two cases: $\xi>0$ and $\xi=0$.

{\it Case 1: $\xi>0$}. Set $d=p(p+1)/2$.  Then $\bY_1, \ldots, \bY_n$ are i.i.d. d-dimensional
random vectors with $E(\bY_1)=\mathbf{0}$ and covariance matrix $\mathbf{V}=\mbox{Cov}(\bY_1)$, where $\mathbf{V}$ is a diagonal matrix with the first $p$ diagonal entries being $\xi$ and the rest $d-p$ diagonal entries being $1$. Set
\[
\bS_n=(s_{11},\ldots,s_{pp}, s_{ij}; 1\le i<j\le
p)^T=\frac{1}{\sqrt{n}}\sum^n_{k=1}\bY_k,
\]
where
\begin{equation}\label{wtos}
s_{ii}=\frac{1}{\sqrt{n}}\sum_{k=1}^n(x_{ki}^2-1)=\frac{1}{\sqrt{n}}(w_{ii}-n),
~~s_{ij}=\frac{1}{\sqrt{n}}\sum^n_{k=1}x_{ki}x_{kj}=\frac{1}{\sqrt{n}}w_{ij}
\end{equation}
for $1\le i\ne j\le p$. Recall the symmetric matrix $\bZ=(\bz_{ij})_{1\le i,j\le p}$, where $\{\bz_{ij}, 1\le i\leq j\le p\}$ are
independent normal random variables with mean zero and
$\Var(\bz_{ij})=1$ for $1\le i<j\le p$ and $\Var(\bz_{ii})=\xi>0$
for $1\le i\le p$. Now we define
\[
\bT=\big(\bz_{11},\ldots,\bz_{pp}, \bz_{ij}; 1\le i<j\le p\big)^T.
\]
To apply Lemma \ref{normalapp} to $\bS_n=n^{-1/2}\sum^n_{k=1}\bY_k,$ the quantity $E(|\bY_1|_2^3)$ is needed to be estimated.

We claim $E(|\bY_1|_2^3)=O(p^3)$. In fact, notice
\[
|\bY_1|_2^2=\sum_{1\le i<j\le
p}x_{1i}^2x_{1j}^2+\sum^p_{i=1}\left(x_{1i}^2-1\right)^2\le
\left(\sum^p_{i=1}x_{1i}^2\right)^2+\sum^p_{i=1}\left(x_{1i}^2-1\right)^2,
\]
 we obtain
\begin{eqnarray*}
E\left(|\bY_1|_2^3\right)
&\le &
\sqrt{2}\left\{E\left(\sum^p_{i=1}x_{1i}^2\right)^3+E\left(\sum^p_{i=1}(x_{1i}^2-1)^2\right)^{3/2}\right\}\\
& \le & C\left\{E\left|\sum^p_{i=1}(x_{1i}^2-1)\right|^3+E\left|\sum^p_{i=1}\left[(x_{1i}^2-1)^2-\xi\right]\right|^{3/2} + p^3+(\xi p)^{3/2}\right\}
\end{eqnarray*}
by using the convexity of function $x^{3/2}$ twice and the notation  $\xi=\mbox{Var}(x_{11}^2)=E(x_{11}^2-1)^2$, where $C$ is a numerical constant free of $n$ and $p$.  Review the assumption that $E(x_{11}^6)<\infty$. Observe that both sums above are sums of i.i.d. random variables with mean $0$ and finite third moment, respectively. The Marcinkiewicz-Zygmund inequality (see, e.g., p. 386
and p. 387 from \cite{CT1997}) guarantees that the first and second expectations above are of order $p^{3/2}$ and $p^{3/4}$, respectively. Hence the claim is true.

Then it follows from Lemma~\ref{normalapp} and $d=p(p+1)/2$  that
\begin{equation}\label{last-onemile}
\sup_{A\in \cal C}|P\left(\bS_n\in A\right)-P\left(\bT\in
A\right)|=O\left(\frac{p^{7/2}}{n^{1/2}}\right)=o(1)
\end{equation}
as $n\to\infty$ by using assumption $p=o(n^{1/7})$.  Define
\[\varphi(x,y,z)=\sqrt{x^2+\frac{1}{4}(y-z)^2}+\frac1{2}(y+z),~~ (x,y,z)\in \mathbb{R}^3,
\]
which is the largest eigenvalue of the symmetric matrix
$$\left(
          \begin{array}{cc}
            y & x \\
            x & z \\
          \end{array}
        \right).
$$
We claim
\begin{equation}\label{convex}
\varphi(x,y,z) \hbox{ is a convex function in } \mathbb{R}^3.
\end{equation}
To check this, it suffices to show $C_t:=\{(x,y,z): \varphi(x,y,z)\le t\} \hbox{ is a convex set in }
\mathbb{R}^3$ for every $t$. In fact, assume $(x,y,x)\in C_t$ and $(x',y',x')\in C_t$ and
$s\in (0,1)$. Set $A=(sx, (1/2)s(y-z))$, $B=((s-1)x', (1/2)(s-1)(y'-z'))$ and $O=(0, 0)$. They are three points in $\mathbb{R}^2$. Use $\|\overrightarrow{AB}\| \leq \|\overrightarrow{OA}\|+\|\overrightarrow{OB}\|$ to see
\begin{eqnarray*}
&&    \sqrt{\big(sx+(1-s)x'\big)^2+\frac14\big(s(y-z)+(1-s)(y'-z')\big)^2}\\
&\le& \sqrt{(sx)^2+\frac14\big(s(y-z)\big)^2}+\sqrt{\big((1-s)x'\big)^2+\frac14\big((1-s)(y'-z')\big)^2}.
\end{eqnarray*}
Therefore,
\begin{eqnarray*}
&&\varphi(s(x,y,z)+(1-s)(x',y',z'))\\
&=&\sqrt{\big(sx+(1-s)x'\big)^2+\frac14\big(s(y-z)+(1-s)(y'-z')\big)^2}\\
&&~~~~+\frac{1}{2}\big[s(y+z)+(1-s)(y'+z')\big]\\
&\le& \sqrt{(sx)^2+\frac14\big(s(y-z)\big)^2}+\sqrt{\big((1-s)x'\big)^2+\frac14\big((1-s)(y'-z')\big)^2}\\
&&~~~+\frac{1}{2}s(y+z)+\frac12(1-s)(y'+z')\\
&=&s\left[\sqrt{x^2+\frac14(y-z)^2}+\frac{1}{2}(y+z)\right]+(1-s)\left[\sqrt{x'^2+\frac14(y'-z')^2}+\frac12(y'+z')\right]\\
&=&s\varphi(x,y,z)+(1-s)\varphi(x',y',z')\\
&\le t.&
\end{eqnarray*}
This implies
$s(x,y,z)+(1-s)(x',y',z')\in C_t.$  We then complete the proof of
\eqref{convex}.

If $x,y,z$ are three coordinates of a vector in $\mathbb{R}^d$ for $d>3$,  we
can treat $\varphi(x,y,x)$ as a special function defined in $\mathbb{R}^d$.
As a consequence, by \eqref{convex}, $L_{ij}:=\varphi(s_{ij},s_{ii},s_{jj})$ is a convex function defined on $\mathbb{R}^d$ with $d=p(p+1)/2$ for any
$1\le i< j\le p$. So is the function $\max_{1\le i<j\le p}L_{ij}$.
Recall that
$l_{ij}=\sqrt{\bz_{ij}^2+\frac{1}{4}(\bz_{ii}-\bz_{jj})^2}+
\frac{1}{2}(\bz_{ii}+\bz_{jj})$, defined in \eqref{Dimsum},
is the same as $\varphi(\bz_{ij},\bz_{ii},\bz_{jj})$, where
$\{\bz_{11},\ldots,\bz_{pp}, \bz_{ij}; 1\le i<j\le p\}$ are the
coordinates of $\bT$. Then it follows from \eqref{last-onemile} that
\[
\lim\limits_{n\to\infty}\sup_{t\in R}\Big|P\Big(\max\limits_{1\le i<j\le
p}L_{ij}\le t\Big)-P\big(\max\limits_{1\le i<j\le p}l_{ij}\le
t\Big)\Big|=0.
\]
For any $z\in \mathbb{R}$, by setting $t=A_{p,
\xi}^{-1}z+B_{p,\xi}$ and applying Theorem~\ref{main1} we have
\begin{equation}\label{almostdone}
\lim\limits_{n\to\infty}P\Big(\max\limits_{1\le i<j\le p}L_{ij}\le B_{p,\xi}+A_{p,
xi}^{-1}z\Big)=G_{\xi}(z).
\end{equation}
Using \eqref{wtos} we have that, for $1\le i<j\le p$,
\begin{eqnarray*}
L_{ij}&=&\varphi\left(\frac{1}{\sqrt{n}}w_{ij},\frac{1}{\sqrt{n}}(w_{ii}-n),
\frac{1}{\sqrt{n}}(w_{jj}-n)\right)\\
&=&\frac{1}{\sqrt{n}}\left(\sqrt{w_{ij}^2+\frac{1}{4}(w_{ii}-w_{jj})^2}+\frac{1}{2}(w_{ii}+w_{jj})-n\right).
\end{eqnarray*}
Since
$[w_{ij}^2+\frac{1}{4}(w_{ii}-w_{jj})^2]^{1/2}+\frac{1}{2}(w_{ii}+w_{jj})$ is
the largest eigenvalue of
\[
 W_{ij}= \left(
          \begin{array}{cc}
            w_{ii} & w_{ij} \\
            w_{ij} & w_{jj} \\
          \end{array}
          \right),
\]
we have from \eqref{NEWS} that
 \[
\frac{T_{n,p,2}-n}{\sqrt{n}}=\max\limits_{1\le i<j\le p}L_{ij}.
 \]
This together with \eqref{almostdone} conclude that
\[
A_{p,\xi}\left(\frac{T_{n,p,2}-n}{\sqrt{n}}-B_{p,\xi}\right)
\]
converges weakly to a probability distribution with distribution  function $G_{\xi}(z)$. The proof is completed.

{\it Case 2: $\xi = 0$.} By assumption, $E(x_{11}^2-1)=0$ and  $\xi=\mbox{Var}(x_{11}^2) =0.$ Thus, $x_{ij}^2=x_{11}^2=1$ almost surely for any $1\leq i, j\leq p$. The formulas from \eqref{NEWS} and \eqref{mooncake} indicate that
$T_{n,p,2}=n+\max_{1\leq i<j \leq} |w_{ij}|.$ Recall \eqref{wtos} and  Theorem~\ref{main1} with $\xi=0$. In {\it Case 1}, change ``$(s_{11},\ldots,s_{pp}, s_{ij}; 1\le i<j\le
p)$" to ``$(s_{ij}; 1\le i<j\le
p)$" and ``$(\bz_{11},\ldots,\bz_{pp}, \bz_{ij}; 1\le i<j\le p\big)$" to ``$(\bz_{ij}; 1\le i<j\le p\big)$", respectively. The same argument in  {\it Case 1} is still valid in the current case. The desired result then follows.  \eop


\section{Proofs of technical  lemmas}\label{lemmas}

In this section, we will prove Lemmas~\ref{lem0}, \ref{rep0}, \ref{rep2}, \ref{lem2} and \ref{estofq}.

\medskip




\noindent{\it Proof of Lemma~\ref{lem0}.} Set $C_p=\{\max_{1\le i\le
p}\bz_i\le \sqrt{2\ln
p}\}=\bigcap^p_{i=1}\{\bz_i\le\sqrt{2\ln p}\}$. Then,
$P(C_p)=\Phi(\sqrt{2\ln p})^p=[1-[1+o(1)](p\sqrt{4\pi\ln
p})^{-1}]^p\to 1$ as $p\to\infty$ by using \eqref{GaussianTail}. Denote
the complement of $C_p$ as $\bar{C}_p$. Then  $P(\bar C_p)\to 0$ as $p\to\infty$. Notice that, on $C_p$,
\[
t_p-\frac{\sqrt{\xi}}{2}(\bz_i+\bz_j)\ge t_p-\sqrt{\xi}\max_{1\le
i\le p}\bz_i\ge t_p-\sqrt{\xi}\sqrt{2\ln p}=\big(r-\sqrt{\xi}+o(1)\big)\sqrt{2\ln
p}\to\infty\]
by using assumption \eqref{tp} and
\begin{eqnarray}
(t_p-\sqrt{\xi}\bz_i)(t_p-\sqrt{\xi}\bz_j)&\ge&
(t_p-\sqrt{\xi}\max\limits_{1\le j\le p}\bz_j)^2\nonumber\\
&\ge&
(t_p-\sqrt{\xi}\sqrt{2\ln p})^2\nonumber\\
&=&2\big(r-\sqrt{\xi}+o(1)\big)^2\ln p\to \infty\label{lowerbound}
\end{eqnarray}
 as $p\to\infty$.

Since $l_{ij}=\sqrt{\bz_{ij}^2+\frac{\xi}{4}(\bz_i-\bz_j)^2}+
\frac{\sqrt{\xi}}{2}(\bz_i+\bz_j)$ for all $1\le i< j\le p$,
inequality $l_{ij}\le t_p$ is equivalent to
\[
\sqrt{\bz_{ij}^2+\frac{\xi}{4}(\bz_i-\bz_j)^2}\le t_p-
\frac{\sqrt{\xi}}{2}(\bz_i+\bz_j).
\]
Since the right-hand side above is non-negative on $C_p$, by squaring on both sides above, we obtain after
simplification that
\[
\{l_{ij}\le t_p\}=\{\bz_{ij}^2\le
(t_p-\sqrt{\xi}\bz_i)(t_p-\sqrt{\xi}\bz_j\}
\]
on $C_p,$ that is,
\begin{equation}\label{onCp}
C_p\bigcap\left\{\max_{1\le i<j\le p}l_{ij}\le
t_p\right\}=C_p\bigcap\limits_{1\le i<j\le p}\left\{\bz_{ij}^2\le
(t_p-\sqrt{\xi}\bz_i)(t_p-\sqrt{\xi}\bz_j\right\}.
\end{equation}
Note that
\begin{eqnarray*}
&&\bigcap\limits_{1\le i<j\le p}\Big\{\bz_{ij}^2\le
(t_p-\sqrt{\xi}\bz_i)(t_p-\sqrt{\xi}\bz_j)\Big\}\\
&=& \Big\{\sum_{1\le
i<j\le p}I\big(\bz_{ij}^2>
(t_p-\sqrt{\xi}\bz_i)(t_p-\sqrt{\xi}\bz_j)\big)=0\Big\}.
\end{eqnarray*}
Recall \eqref{q(x,y)}. Conditional on $\bz_1, \bz_2, \ldots, \bz_p$, $\{\bz_{ij}^2>
(t_p-\sqrt{\xi}\bz_i)(t_p-\sqrt{\xi}\bz_j)\}$, $1\le i<j\le p$, are
$p(p-1)/2$ independent events with conditional probabilities
$q(\bz_i,\bz_j;t_p)$. Write
\begin{eqnarray*}
\Delta_p
&=&P\Big(\bigcap\limits_{1\le i<j\le p}\{\bz_{ij}^2\le
(t_p-\sqrt{\xi}\bz_i)(t_p-\sqrt{\xi}\bz_j)\}\big|\bz_1,\bz_2,
\ldots, \bz_p\Big)\\
&& -\exp\Big(-\sum_{1\le i<j\le p}q(\bz_i,\bz_j;t_p)\Big).
\end{eqnarray*}
By using the Poisson approximation
\begin{eqnarray*}
|\Delta_p|&\le& \min\Big\{1, \Big(\sum_{1\le i<j\le
p}q(\bz_i,\bz_j;t_p)\Big)^{-1}\Big\}\sum_{1\le i<j\le
p}q(\bz_i,\bz_j;t_p)^2\\
&\le& \max_{1\le i<j\le p}q(\bz_i,\bz_j;t_p).
\end{eqnarray*}
See, e.g., page 8 in Barbour {\it et al.}~\cite{BHJ1992}.  Then  we have from \eqref{GaussianTail}, \eqref{q(x,y)} and \eqref{lowerbound} that
\[
\max_{1\le i<j\le p}q(\bz_i,\bz_j;t_p)\le
2\big[1-\Phi(t_p-\sqrt{\xi}\sqrt{2\ln p})\big]
\to 0
\]
on $C_p$, which yields $E|\Delta_p|\to 0$ as $p\to\infty$. Therefore,  it follows
from \eqref{onCp} that
\begin{eqnarray*}
    &&P\big(\max_{1\le i<j\le p}l_{ij}\le t_p\big)\\
&=&
P\Big(C_p\bigcap\big\{\max_{1\le i<j\le p}l_{ij}\le t_p\big\}\Big)+P\big(\bar
C_p\bigcap\big\{\max_{1\le i<j\le
p}l_{ij}\le t_p\big\}\big)\\
&=&P\Big(C_p\bigcap\limits_{1\le i<j\le p}\{\bz_{ij}^2\le
(t_p-\sqrt{\xi}\bz_i)(t_p-\sqrt{\xi}\bz_j)\}\Big)+O(P(\bar C_p))\\
&=&EP\Big(C_p\bigcap\limits_{1\le i<j\le p}\{\bz_{ij}^2\le
(t_p-\sqrt{\xi}\bz_i)(t_p-\sqrt{\xi}\bz_j)\}\big|\bz_1,\bz_2,
\ldots, \bz_p\Big)+o(1)\\
&=&E\Big[I(C_p)P\Big(\bigcap\limits_{1\le i<j\le p}\{\bz_{ij}^2\le
(t_p-\sqrt{\xi}\bz_i)(t_p-\sqrt{\xi}\bz_j)\}\big|\bz_1,\bz_2,
\ldots, \bz_p\Big)\Big]+o(1).
\end{eqnarray*}
By the estimate of $\Delta_p$, the above is identical to
\begin{eqnarray*}
&&E\Big[I(C_p)\exp\Big(-\sum_{1\le i<j\le
p}q(\bz_i,\bz_j;t_p)\Big)+\Delta_p\Big]+o(1)\\
&=&E\exp\Big(-\sum_{1\le i<j\le p}q(\bz_i,\bz_j;t_p)\Big)-
E\Big[I(\bar C_p)\exp\Big(-\sum_{1\le i<j\le
p}q(\bz_i,\bz_j;t_p)\Big)\Big]\\
&&+E\big[I(C_p)\Delta_p\big]+o(1)\\
&=&E\exp\Big(-\sum_{1\le i<j\le p}q(\bz_i,\bz_j;t_p)\Big)+
O(P(\bar C_p))+ o(1)\\
&=&E\exp\Big(-\sum_{1\le i<j\le p}q(\bz_i,\bz_j;t_p)\Big)+ o(1).
\end{eqnarray*}
This completes the proof of  Lemma~\ref{lem0}. \eop

\vspace{10pt}


\noindent{\it Proof of Lemma~\ref{rep0}}. For each $1\le i\le p$, define
$M_i=\{\max_{1\le j\le p, \,j\ne
i}\bz_j<\bz_i\big\}$.  By using the symmetry, we have
\begin{eqnarray*}
&&Ef_p(\bz_1, \cdots, \bz_{p-1}, \bz_p)\\
&=&\sum^p_{i=1}E\big[f_p(\bz_1, \cdots, \bz_{p-1},
\bz_p)I(M_i)\big]\\
&=&p\cdot E\big[f_p(\bz_1, \cdots, \bz_{p-1}, \bz_p)I(M_p)\big]\\
&=&p\cdot E\left\{\int^{\bz_{p}}_{-\infty}\cdots\int^{\bz_{p}}_{-\infty} f_p(z_1, \cdots, z_{p-1}, \bz_p)\Big[\prod^{p-1}_{j=1}\phi(z_j)\Big]dz_1\cdots dz_{p-1}\right\}\\
&=&p\int^\infty_{-\infty}\phi(x)\left\{\int^x_{-\infty}\cdots\int^x_{-\infty} f_p(z_1, \cdots, z_{p-1}, x)\Big[\prod^{p-1}_{j=1}\phi(z_j)\Big]dz_1\cdots dz_{p-1}\right\}dx.
\end{eqnarray*}
Use \eqref{tongzhia} to rewrite the above as
\begin{eqnarray*}
&&p\int^\infty_{-\infty}\phi(x)\Big[\int_{\mathbb{R}^{p-1}}f_p(z_1, \cdots, z_{p-1}, x)\Big(\prod^{p-1}_{j=1}\phi(z_j)I(z_j\leq x)\Big)dz_1\cdots dz_{p-1}\Big]dx\\
&=& p\int^\infty_{-\infty}\Phi(x)^{p-1}\phi(x)\Big[\int_{\mathbb{R}^{p-1}}f_p(z_1, \cdots, z_{p-1}, x)\prod^{p-1}_{j=1}\phi_x(z_j)dz_1\cdots dz_{p-1}\Big]dx\\
&=&\int^\infty_{-\infty}\lambda_p(y)\Big[\int_{\mathbb{R}^{p-1}}
f_p(z_1, \cdots, z_{p-1},
c_p(y))\prod^{p-1}_{j=1}\phi_{c_p(y)}(z_j)dz_1\cdots dz_{p-1}\Big]
dy
\end{eqnarray*}
by setting $x=c_p(y)=\beta_p+\alpha_p^{-1}y$ in the last step,
where
$\lambda_p(y):=p\alpha_p^{-1}\phi(c_p(y))\Phi(c_p(y))^{p-1}.$ It is easy to see that  $\lambda_p(y)$
is the density function of $\alpha_p(\max_{1\le j\le
p}\bz_j-\beta_p)$.  Observe that the integral  within the square bracket above is equal to $Ef_p(\bar\bz_1, \cdots,
\bar\bz_{p-1}, c_p(y))$,  where for each $y\in \mathbb{R}$, $\bar\bz_1, \cdots,
\bar\bz_{p-1}$ are i.i.d. random variables with density function
$\phi_{c_p(y)}(\cdot).$ Therefore,
\begin{equation}\label{nansile}
Ef_p(\bz_1, \cdots, \bz_{p-1}, \bz_p)
=\int^\infty_{-\infty}\lambda_p(y)\cdot Ef_p(\bar\bz_1, \cdots, \bar\bz_{p-1},
c_p(y))dy.
\end{equation}
By \eqref{standard}, $\beta_p=\alpha_p-\alpha_p^{-1}\ln\big(\sqrt{2\pi}\alpha_p\big)~~\mbox{
and}~~~c_p(y)=\beta_p+\alpha_p^{-1}y$. Thus, $\alpha_p/\beta_p\to 1$. Now, by using the notation $\alpha_p=\sqrt{2\ln p}$, it is not difficult to check
\[
p\alpha_p^{-1}\phi(c_p(y))=\frac{p}{\sqrt{2\pi}\,\alpha_p}\cdot \exp\left(-\frac{1}{2}\left(\beta_p+\alpha_p^{-1}y\right)^2\right)\to e^{-y}.
\]
This and \eqref{EV} yield that $\lambda_p(y) \to \lambda(y)=\exp(-y-e^{-y})$ for every $y\in \mathbb{R}$. By the Scheff\'e lemma (see, e.g., p. 55 from \cite{DW1991}), we have
\[
\lim\limits_{p\to\infty}\int^\infty_{-\infty}|\lambda_p(y)-\lambda(y)|\,dy=0.
\]
Consequently,  we have from \eqref{nansile} that
\begin{eqnarray*}
&&\Big|Ef_p\left(\bz_1, \cdots, \bz_{p-1},
\bz_p\right)-\int^\infty_{-\infty}Ef_p\left(\bar\bz_1, \cdots,
\bar\bz_{p-1}, c_p(y)\right)\lambda(y)dy\Big|\\
&=& \Big|\int^\infty_{-\infty}Ef_p\left(\bar\bz_1, \cdots,
\bar\bz_{p-1}, c_p(y)\right)\big[\lambda_p(y)-\lambda(y)\big]dy\Big|\\
&\le& b\int^\infty_{-\infty}|\lambda_p(y)-\lambda(y)|dy\\
&\to& 0
\end{eqnarray*}
as $p\to\infty$. The proof is completed.   \eop

\medskip

\medskip

\noindent{\it Proof of Lemma~\ref{rep2}}. Recall \eqref{QU}. We know
\begin{equation*}
Q_p=\sum_{1\le i<j\le p-1}q(\bar\bz_i, \bar\bz_j;t_p)~~ \mbox{and}~~
U_p=\sum^{p-1}_{j=1}q(c_p, \bar\bz_j;t_p)
\end{equation*}
where $\bar\bz_1, \cdots,
\bar\bz_{p-1}$ are i.i.d. random variables with density function
$\phi_{c_p(y)}(\cdot)$  with $c_p=c_p(y)=\beta_p+\alpha_p^{-1}y$ and
\begin{equation}\label{dayu}
q(\bar\bz_i, \bar\bz_j;t_p) = 2\Big[1-\Phi\Big(\sqrt{(t_p-\sqrt{\xi}\bar\bz_i)(t_p-\sqrt{\xi}\bar\bz_j)}\,\Big)\Big]
\end{equation}
provided the term under the square-root sign is non-negative; see the explanation under \eqref{q(x,y)}. By definition,  $\bar\bz_i \leq c_p(y)$ for each $i.$ Also notice \eqref{tp} is equivalent to that
\begin{equation}\label{gongming}
\lim\limits_{p\to\infty}\frac{t_p}{c_p}=r\in \left(\sqrt{\xi},
\sqrt{\xi}+\frac{2}{\sqrt{\xi}}\right).
\end{equation}
We see that, as $p$ is sufficiently large,  $t_p-\sqrt{\xi}\bar\bz_i\geq t_p-\sqrt{\xi}c_p(y)\geq 0$ a.s. for each $i$. Thus, for the discussion in the future, we always assume \eqref{dayu} holds.

(a). Note that $Q_p$ is a U-statistic and
$U_p$ is a sum of $p-1$ i.i.d. random variables. For $Q_p$, we perform  the classical H\'ajek's projection for $U$-statistic next (see,
e.g. Denker~\cite{Denker85}). To do so, for $1\le i<j\le p-1$, define
\[
v_{ij}=q(\bar\bz_i,\bar\bz_j;t_p)-q(\bar\bz_i;t_p)-q(\bar\bz_j;t_p)+q(t_p),
\]
where $q(x;t_p)=Eq(x,\bar\bz_2;t_p)$ for $x\le c_p$ as defined in \eqref{q(x)} and $q(t_p)=Eq(\bar\bz_1,\bar\bz_2; t_p)=
Eq(\bar\bz_1;t_p)$ in \eqref{qtp}. It is easy to check that  $\{v_{ij}:~1\le i<j\le p\}$ are uncorrected random variables with mean zero and $E(v_{ij}^2)=\Var(q(\bar\bz_1,\bar\bz_1;t_p))$. Thus,
\[ E(v_{ij}^2)\le
E\big[q(\bar\bz_1,\bar\bz_2;t_p)^2\big]\le q(c_p,c_p;t_p)q(t_p)
\]
for any $i,j$, where the inequality $q(\bar\bz_1,\bar\bz_2;t_p)\le q(c_p,c_p;t_p)~ a.s.$ from \eqref{dayu} is used in the last step. Furthermore, to handle $U_p$, we set
\[
v_{ip}=q(c_p,\bar\bz_i;t_p)-q(c_p;t_p)
\]
 for $1\le i<p$. By the same argument used in the estimate of $E(v_{ij}^2)$, we have
\[
\Var(v_{ip})= \Var(q(c_p,\bar\bz_i;t_p))\le
E\big[q(c_p,\bar\bz_i;t_p)^2\big]\le q(c_p,c_p;t_p)q(c_p;t_p).
\]
Define
\[
R_p=\sum\limits_{1\le i<j\le
p-1}v_{ij}+\sum\limits^{p-1}_{i=1}v_{ip}.
\]
Note that $q(x,y;t_p)$ given in \eqref{q(x,y;tp)} is strictly
increasing in both $x$ and $y$ when $x\le c_p$ and $y\le c_p$. Under
conditions \eqref{tp} and \eqref{cp}, there exists $\delta_1>0$ such that
\begin{equation*}
q(c_p,c_p; t_p)\le 2\big[1-\Phi(\delta_1c_p)\big]\to 0~\mbox{ and }q(t_p)\le
2\big[1-\Phi(\delta_1c_p)\big]\to 0
\end{equation*}
as $p\to\infty$. Easily, $\mbox{Cov}(v_{ij}, v_{kp})=0$ for any $1\leq i<j \leq p-1$ and $1\leq k \leq p-1$. Hence,
\[
E\big(R_p^2\big)=\sum_{1\le i<j\le
p-1}\Var(v_{ij})+\sum^{p-1}_{i=1}\Var(v_{ip})=o\Big(p^2q(t_p)+pq(c_p;t_p)\Big)=o(1)
\]
under assumptions that both $pq(c_p;t_p)$ and $p^2q(t_p)$ are
bounded.  Also notice that
\[
R_p=Q_p+U_p-\Big[(p-2)\sum^{p-1}_{j=1}q(\bar\bz_j;t_p)-d_p\Big]
\]
where $d_p=(1/2)(p-1)(p-2)q(t_p)-(p-1)q(c_p;t_p)$. Observe $Q_p+U_p\ge 0$ and
$(p-2)\sum^{p-1}_{j=1}q(\bar\bz_j;t_p)-d_p\ge -d_p\ge c$ for some
$c$ since $d_p$ is bounded.  By using inequality
\[
|e^{-s}-e^{-t}|=\Big|\int^s_te^{-x}dx\Big|\le e^{-a}|s-t|
\]
for $s,t\ge a$, we obtain
\begin{eqnarray*}
&&\Big|Ee^{-(Q_p+U_p)}-E\exp\Big(d_p-(p-2)\sum^{p-1}_{j=1}q(\bar\bz_j;t_p)\Big)\Big|\\
&\le&
e^{|c|}E\Big|Q_p+U_p+d_p-(p-2)\sum^{p-1}_{j=1}q(\bar\bz_j;t_p)\Big|\\
&=&e^{|c|}E|R_p| \\
&\le& e^{|c|} [E(R_p^2)]^{1/2}\\
&\to& 0
\end{eqnarray*}
as $p\to\infty$, proving \eqref{rep2a}.

(b).   We first cite the Taylor reminder theorem. Fix a positive
integer $k$.  For every $x>0$, there exists $\omega\in [0, x]$
such that
\[
e^{-x}=\sum^k_{j=0}\frac{(-x)^j}{j!}+e^{-\omega}\frac{(-x)^{k+1}}{(k+1)!}
\]
which implies
\begin{equation}\label{taylor}
\left|e^{-x}-\sum^k_{j=0}(-1)^j\frac{x^j}{j!}\right|\le
\frac{x^{k+1}}{(k+1)!},~~~x\ge 0,
\end{equation}
for any integer $k\ge 1$. Plug $x=(p-2)q(\bar\bz_1;t_p))$ in \eqref{taylor} for $k=1$ and
then take expectations on both sides to get
\begin{eqnarray*}
&&\big|E\exp\big(-(p-2)q(\bar\bz_1;t_p)\big)-1+
(p-2)Eq(\bar\bz_1;t_p)\big|\\
& \leq & (p-2)^2 E\big[q(\bar\bz_1;t_p)^2\big]\\
&\le &  (p-2)q(c_p;t_p)\cdot
(p-2)Eq(\bar\bz_1;t_p),
\end{eqnarray*}
i.e.,
\[
\left|E\exp\big(-(p-2)q(\bar\bz_1;t_p)\big)-1+ (p-2)q(t_p)\right|\le
(p-2)q(c_p;t_p)\cdot(p-2)q(t_p).
\]
Using the given conditions in part (b) we have $d_p=[1+o(1)]\gamma$ and
\[
E\exp\big(-(p-2)q(\bar\bz_1;t_p)\big)=1-\frac{2[1+o(1)]\gamma}{p}
\]
as $p\to\infty$. Therefore, by using the independence we get
\begin{eqnarray*}
E\exp\left(d_p-(p-2)\sum^{p-1}_{j=1}q(\bar\bz_j;t_p)\right)
&=&e^{d_p}\Big[E\exp\big(-(p-2)q(\bar\bz_1;t_p)
\big) \Big]^{p-1}\\
&=& e^{[1+o(1)]\gamma}\left[1-\frac{2[1+o(1)]\gamma}{p} \right]^{p-1}\\
&\to& e^{-\gamma}.
\end{eqnarray*}
Then \eqref{rep2b} follows from \eqref{rep2a}.

(c). The right-hand side of \eqref{irregularlimit} is well defined
under given conditions. The series
\[
\kappa:=\sum^\infty_{j=1}(-1)^{j-1}\frac{\gamma_j}{j!}
\]
is absolutely convergent since
\[
\gamma_j\le
\lim\limits_{p\to\infty}p^2\big[pq(c_p;t_p)\big]^{j-1}Eq(\bar\bz_1;t_p)=
\lim\limits_{p\to\infty}\big[pq(c_p;t_p)\big]^{j-1}p^2q(t_p)=2\tau^{j-1}\gamma_1.
\]
Given $k\ge 2$, we apply \eqref{taylor} with
$x=(p-2)\bar\bz_1$ and take expectations on both sides to get
\begin{eqnarray*}
&& \left|E\exp\big(-(p-2)q(\bar\bz_1;t_p)\big)-1-\sum^k_{j=1}
\frac{(-1)^j(p-2)^jE\big(q(\bar\bz_1;t_p)^j\big)}{j!}\right|\\
&\le &
\frac{(p-2)^{k+1}E\big(q(\bar\bz_1;t_p)^{k+1}\big)}{(k+1)!}.
\end{eqnarray*}
By using assumptions $p^{j+1}E\big(q(\bar\bz_1;t_p)^j\big)\to
\gamma_j$ for $j\ge 2$ and $p^2Eq(\bar\bz_1;t_p)\to
2\gamma_1$ we have
\[
\limsup_{p\to\infty}\left|p\Big[E\exp\big(-(p-2)q(\bar\bz_1;t_p)\big)-1\Big]+2\gamma_1+\sum^k_{j=2}
\frac{(-1)^{j-1}\gamma_j}{j!}\right|\le \frac{\gamma_{k+1}}{(k+1)!}\le
\frac{2\tau^k\gamma_1}{(k+1)!}.
\]
Thus, we have
\begin{eqnarray*}
&&\limsup_{p\to\infty}\left|p\Big[E\exp\big(-(p-2)q(\bar\bz_1;t_p)\big)-1\Big]+\gamma_1+\kappa\right|\\
&\le&
\limsup_{p\to\infty}\left|p\Big[E\exp\big(-(p-2)q(\bar\bz_1;t_p)\big)-1\Big]+\gamma_1+\sum^k_{j=1}
\frac{(-1)^{j-1}\gamma_j}{j!}\right|+\left|\sum^{\infty}_{j=k+1}\frac{(-1)^{j-1}\gamma_j}{j!}\right|\\
&\le&
\frac{2\tau^k\gamma_1}{(k+1)!}+\sum^{\infty}_{j=k+1}\frac{\gamma_j}{j!}.
\end{eqnarray*}
By letting $k\to\infty$, we get
\[
\lim\limits_{p\to\infty}\left|p\Big[E\exp(-(p-2)q(\bar\bz_1;t_p)\big)-1\Big]+\gamma_1+\kappa\right|=0,
\]
which implies
\[
\left[E\exp\big(-(p-2)q(\bar\bz_1;t_p) \big) \right]^{p-1}\to
\exp(-\gamma_1-\kappa).
\]
Therefore,
\begin{eqnarray*}
E\exp\left(d_p-(p-2)\sum^{p-1}_{j=1}q(\bar\bz_j;t_p)\right)
&=&e^{d_p}\Big[E\exp(-(p-2)q(\bar\bz_1;t_p)
\big)\Big]^{p-1}\\
&\to& e^{\gamma_1-\tau}\exp(-\gamma_1-\kappa)\\
&=&\exp(-\tau-\kappa).
\end{eqnarray*}
Consequently, we obtain \eqref{irregularlimit} from \eqref{rep2a}.
    \eop

\medskip

\noindent{\it Proof of Lemma~\ref{lem2}}. Recall \eqref{gongming}. We have from
\eqref{tp} and \eqref{cp} that there exist constants $\delta_1, \delta_2>0$ such that
\begin{equation}\label{tp-range}
\delta_1c_p\le t_p-\sqrt{\xi} c_p\le\delta_2c_p.
\end{equation}
Then, for
$x,y\le c_p$, we have  $t_p-\sqrt{\xi}x\ge t_p-\sqrt{\xi}c_p\ge
\delta_1c_p$ and $t_p-\sqrt{\xi}y\ge t_p-\sqrt{\xi}c_p\ge
\delta_1c_p$ if $p$ is large.  The assertion \eqref{GaussianTail} implies
\[
1-\Phi\Big(\sqrt{(t_p-\sqrt{\xi}x)(t_p-\sqrt{\xi}y)}\,\Big)=\frac{[1+o(1)]\exp\big(-\frac12\big(t_p-\sqrt{\xi}x\big)\big(t_p-\sqrt{\xi}y\big)\big)
} {\sqrt{2\pi}\sqrt{\big(t_p-\sqrt{\xi}x\big)\big(t_p-\sqrt{\xi}y\big)}}
\]
uniformly  over $x, y\le c_p$ as $p\to\infty$.  By \eqref{q(x,y;tp)} and \eqref{q(x)}, we see
\begin{eqnarray}
q(x;t_p)
&=&\frac{1}{\Phi(c_p)}\int^{c_p}_{-\infty}2\Big(1-\Phi\Big(\sqrt{(t_p-\sqrt{\xi}x)(t_p-\sqrt{\xi}y)}\Big)\Big)\phi(y)dy\nonumber\\
&=&\frac{2[1+o(1)]}{\Phi(c_p)\sqrt{2\pi}}\int^{c_p}_{-\infty}\frac1{\sqrt{(t_p-\sqrt{\xi}x)(t_p-\sqrt{\xi}y)}}
\exp\left(-\frac12\big(t_p-\sqrt{\xi}x\big)\big(t_p-\sqrt{\xi}y\big)\right)\phi(y)dy\nonumber\\
&=&\frac{2[1+o(1)]}{\sqrt{2\pi}\sqrt{t_p-\sqrt{\xi}x}}\int^{c_p}_{-\infty}\frac1{\sqrt{2\pi}\sqrt{t_p-\sqrt{\xi}y}}
\exp\left(-\frac12\big((t_p-\sqrt{\xi}x)(t_p-\sqrt{\xi}y)+y^2\big)\right)dy\nonumber\\
&=&:\frac{2[1+o(1)]}{\sqrt{2\pi}\sqrt{t_p-\sqrt{\xi}x}}\cdot
I_1\label{q(x)=},
\end{eqnarray}
where $I_1$ denotes the integral in the equation above \eqref{q(x)=}.
To evaluate the integral $I_1$,  we see
that, for any $y\in \mathbb{R}$,
\begin{eqnarray}
&&\big(t_p-\sqrt{\xi}x\big)\big(t_p-\sqrt{\xi}y\big)+y^2\nonumber\\
&=&\Big[y-\frac12\sqrt{\xi}\big(t_p-\sqrt{\xi}x\big)\Big]^2+\big(t_p-\sqrt{\xi}x\big)t_p-\frac{1}{4}\xi\big(t_p-\sqrt{\xi}x\big)^2\nonumber\\
&=&\Big[y-\frac12\sqrt{\xi}\big(t_p-\sqrt{\xi}x\big)\Big]^2+
\frac{4-\xi}{4}t_p^2-\frac{\xi^2}{4}x^2-\frac{(2-\xi)\sqrt{\xi}}{2}xt_p,\label{wecan}
\end{eqnarray}
and then we obtain
\begin{equation}\label{Ip}
I_1=\exp\Big(-\frac12\Big(\frac{4-\xi}{4}t_p^2-\frac{\xi^2}{4}x^2-\frac{(2-\xi)\sqrt{\xi}}{2}xt_p\Big)
\Big)\cdot I_2,
\end{equation}
where
\begin{equation*}
I_2:=\int^{c_p}_{-\infty}\frac{\exp\big(-\frac12\big(y-\frac12\sqrt{\xi}(t_p-\sqrt{\xi}x)\big)^2\big)}{\sqrt{2\pi}\sqrt{t_p-\sqrt{\xi}y}}dy.
\end{equation*}
Define $s$ and $\bar c_p$ such that
$s=y-\frac12\sqrt{\xi}(t_p-\sqrt{\xi}x)=y-(c_p-\bar c_p)$. Evidently, $\bar c_p=c_p-\frac12\sqrt{\xi}(t_p-\sqrt{\xi}x)$. Then after changing variables, we see
\begin{eqnarray}
I_2
&=&\int^{\bar
c_p}_{-\infty}\frac{\exp\big(-s^2/2\big)}{\sqrt{2\pi}\sqrt{t_p-\sqrt{\xi}(s-\bar
c_p+c_p)}}ds\label{integral-zero}\\
&=&\int^{\bar
c_p}_{-\infty}\frac{\exp\big(-s^2/2\big)}{\sqrt{2\pi}\sqrt{t_p-\frac{\xi}{2}(t_p-\sqrt{\xi}x)-\sqrt{\xi}s}}ds\nonumber\\
&=&\int^{\bar
c_p}_{-\infty}\frac{\exp\big(-s^2/2\big)}{\sqrt{2\pi}\sqrt{(1-0.5\xi)t_p+0.5\xi^{3/2}x
-\sqrt{\xi}s}}ds\label{integral-one}.
\end{eqnarray}
Now we start to prove parts (i), (ii) and (iii) stated in the lemma.

{\it (i)}.  It is easy to verify that $(1-0.5\xi)t_p+0.5\xi^{3/2}x
-\sqrt{\xi}s\ge t_p-\sqrt{\xi}c_p$ for every $s\le \bar c_p$ and $x\leq c_p$. We have from
\eqref{integral-one} that
\[
I_2\le\frac{1}{\sqrt{t_p-\sqrt{\xi}c_p}}\int^{\bar
c_p}_{-\infty}\frac{\exp\big(-s^2/2\big)}{\sqrt{2\pi}}ds
\le\frac{\Phi(\bar c_p)}{\sqrt{t_p-\sqrt{\xi}c_p}}
\]
which, coupled with  \eqref{q(x)=} and \eqref{Ip}, implies
\eqref{q(x)-upper}.

{\it (ii)}.  Note that
$\bar c_p=\bar c_p(t_p,x)\ge \bar
c_p(t_p,x_p)=(1/2)\xi b_p$ for $x_p\le x\le c_p$.  We have from \eqref{tp-range} that
$(1-0.5\xi)t_p+0.5\xi^{3/2}x -\sqrt{\xi}s\ge t_p-\sqrt{\xi}c_p \ge
\delta_1c_p\to\infty$ uniformly over $s\le \bar c_p$, and
$(1-0.5\xi)t_p+0.5\xi^{3/2}x\geq t_p-\sqrt{\xi}c_p+ 0.5\xi^{3/2}b_p\ge
\delta_1c_p+ 0.5\xi^{3/2}b_p\to\infty$.
 Therefore, we have
\begin{eqnarray}
&&\left|\frac1{\sqrt{(1-0.5\xi)t_p+0.5\xi^{3/2}x
-\sqrt{\xi}s}}-\frac1{\sqrt{(1-0.5\xi)t_p+0.5\xi^{3/2}x}}\right| \nonumber\\
&=&\frac{\sqrt{\xi}|s|}{\sqrt{(1-0.5\xi)t_p+0.5\xi^{3/2}x
-\sqrt{\xi}s}\sqrt{(1-0.5\xi)t_p+0.5\xi^{3/2}x}}\nonumber\\
&&~~~\cdot\frac1{\sqrt{(1-0.5\xi)t_p+0.5\xi^{3/2}x
-\sqrt{\xi}s}+\sqrt{(1-0.5\xi)t_p+0.5\xi^{3/2}x}}\nonumber\\
&\le&\frac{\sqrt{\xi}
|s|}{2(t_p-\sqrt{\xi}c_p)\sqrt{(1-0.5\xi)t_p+0.5\xi^{3/2}x}}\label{France}
\end{eqnarray}
uniformly for $s\le \bar c_p$ and  $x_p\le x\le c_p$.
Notice
\begin{eqnarray*}
&&\left|I_2 -\frac{\Phi(\bar c_p)}{\sqrt{(1-0.5\xi)t_p+0.5\xi^{3/2}x
}}\right|\\
&=&\left|\int^{\bar
c_p}_{-\infty}\frac{\phi(s)}{\sqrt{(1-0.5\xi)t_p+0.5\xi^{3/2}x
-\sqrt{\xi}s}}ds -\int^{\bar
c_p}_{-\infty}\frac{\phi(s)}{\sqrt{(1-0.5\xi)t_p+0.5\xi^{3/2}x
}}ds\right|\\
&\le& \int^{\bar
c_p}_{-\infty}\left|\frac1{\sqrt{(1-0.5\xi)t_p+0.5\xi^{3/2}x
-\sqrt{\xi}s}}-\frac1{\sqrt{(1-0.5\xi)t_p+0.5\xi^{3/2}x
}}\right|\phi(s)ds.
\end{eqnarray*}
By \eqref{France}, the above is dominated by
\begin{eqnarray*}
& & \int^{\bar c_p}_{-\infty}\frac{\sqrt{\xi}
|s|\phi(s)}{(t_p-\sqrt{\xi}c_p)\sqrt{(1-0.5\xi)t_p+0.5\xi^{3/2}x}}ds\\
&\le &
\frac{\sqrt{\xi}}{(t_p-\sqrt{\xi}c_p)\sqrt{(1-0.5\xi)t_p+0.5\xi^{3/2}x
}} \int^{\infty}_{-\infty}|s|\phi(s)ds \\
&\le
&\frac{\sqrt{\xi}}{(t_p-\sqrt{\xi}c_p)\sqrt{(1-0.5\xi)t_p+0.5\xi^{3/2}x}}
\end{eqnarray*}
where in the last step we use the inequality $\int^{\infty}_{-\infty}|s|\phi(s)ds=E|N(0,1)| \leq [EN(0,1)^2]^{1/2}=1$.
Since $1\ge \Phi(\bar c_p)\ge \Phi(\frac{\xi}{2}b_p)\to 1$ uniformly
as $p\to\infty$, we obtain
\begin{equation}\label{proofofI2}
 I_2=\frac{1+o(1)}{\sqrt{(1-0.5\xi)t_p+0.5\xi^{3/2}x }}
\end{equation}
holds uniformly over $x_p\le x\le c_p$ as $p\to\infty$. Therefore,
from \eqref{q(x)=}, \eqref{Ip}, and \eqref{integral-one} we conclude
\eqref{q(x)-approx}.


{\it (iii)}.  We now evaluate $I_2$ defined in \eqref{integral-zero}. In fact,
\begin{eqnarray*}
I_2&=&\int^{\bar
c_p}_{-\infty}\frac{\exp(-s^2/2)}{\sqrt{2\pi}\sqrt{t_p-\sqrt{\xi}c_p+\sqrt{\xi}(\bar
c_p-s)}}ds\\
&=&
-\int^{\bar
c_p}_{-\infty}
\frac{\big(\exp(-s^2/2)\big)'}{s\sqrt{2\pi}\sqrt{t_p-\sqrt{\xi}c_p+\sqrt{\xi}(\bar
c_p-s)}}ds.
\end{eqnarray*}
By integration by parts, the above is identical to
\begin{eqnarray*}
&&-\frac{\exp\big(-s^2/2\big)}{s\sqrt{2\pi}\sqrt{t_p-\sqrt{\xi}c_p+\sqrt{\xi}(\bar
c_p-s)}}\Big|^{\bar c_p}_{-\infty}\\
&&~~~+\int^{\bar
c_p}_{-\infty}\frac{\exp(-s^2/2)}{\sqrt{2\pi}}\frac{d}{ds}\Big(\frac1{s\sqrt{t_p-\sqrt{\xi}c_p+\sqrt{\xi}(\bar
c_p-s)}}\Big)ds\\
&=&-\frac{\exp(-\bar c_p^2/2)}{\bar
c_p\sqrt{2\pi}\sqrt{t_p-\sqrt{\xi}c_p}}+\frac{\sqrt{\xi}}{2}\int^{\bar
c_p}_{-\infty}\frac{\exp\big(-s^2/2\big)}{\sqrt{2\pi}s\big[t_p-\sqrt{\xi}c_p+\sqrt{\xi}(\bar
c_p-s)\big]^{3/2}}ds\\
&&~~~-\int^{\bar
c_p}_{-\infty}\frac{\exp\big(-s^2/2\big)}{\sqrt{2\pi}s^2\big[t_p-\sqrt{\xi}c_p+\sqrt{\xi}(\bar
c_p-s)\big]^{1/2}}ds\\
&=&:-\frac{\exp\big(-\bar c_p^2/2)}{\bar
c_p\sqrt{2\pi}\sqrt{t_p-\sqrt{\xi}c_p}}  +I_3+I_4.
\end{eqnarray*}
Note that
$\bar c_p=
\bar c_p(t_p,x)=c_p-\frac12\sqrt{\xi}(t_p-\sqrt{\xi}x)$ is increasing in $x$, we have
 $\bar c_p\le \bar
c_p(t_p,x_p)= -\frac{\xi}{2}b_p\to-\infty$ uniformly over $x\le
\bar x_:=\frac1{\sqrt{\xi}}\big(t_p-\frac{2c_p}{\sqrt{\xi}}\big)-b_p$ as
$p\to\infty$. If $s\le \bar c_p$, we have $|s|\ge -\bar
c_p\ge \frac{\xi}{2}b_p$,
\[
t_p-\sqrt{\xi}c_p+\sqrt{\xi}(\bar c_p-s)\ge
t_p-\sqrt{\xi}c_p\ge\delta_1c_p\to\infty,
\]
and we can easily show that $I_3=o(I_2)$ and  $I_4=o(I_2)$.  Then we
conclude
\begin{equation}\label{I2}
I_2=-[1+o(1)]\frac{\exp\big(-\frac12\bar c_p^2\big)}{\sqrt{2\pi}\bar
c_p\sqrt{t_p-\sqrt{\xi}c_p}}.
\end{equation}
Now we have from a straightforward calculation that
\begin{eqnarray*}
J:&=&-\frac{\exp\big(-\frac12\bar c_p^2\big)}{\bar
c_p}\exp\Big(-\frac12\Big(\frac{4-\xi}{4}t_p^2-\frac{\xi^2}{4}x^2-\frac{(2-\xi)\sqrt{\xi}}{2}xt_p\Big)\Big)\\
&=&\frac{\exp(-\frac12c_p^2)}{\frac12\sqrt{\xi}\big(t_p-\sqrt{\xi}x\big)-c_p}\exp\Big(-\frac{1}{2}\big(t_p-\sqrt{\xi}x\big)\big(t_p-\sqrt{\xi}c_p\big)\Big).
\end{eqnarray*}
This can also be obtained from \eqref{wecan} by setting $y=c_p$ and using  $\bar c_p=c_p-\frac12\sqrt{\xi}(t_p-\sqrt{\xi}x)$. It follows from \eqref{q(x)=}, \eqref{Ip} and \eqref{I2}
that
\begin{eqnarray*}
q(x;t_p)&=&\frac{2[1+o(1)]}{\sqrt{2\pi}\sqrt{t_p-\sqrt{\xi}x}}\cdot
\exp\Big(-\frac12\Big(\frac{4-\xi}{4}t_p^2-\frac{\xi^2}{4}x^2-\frac{(2-\xi)\sqrt{\xi}}{2}xt_p\Big)
\Big)\cdot I_2\\
&=&\frac{1+o(1)}{\pi\sqrt{t_p-\sqrt{\xi}c_p}\sqrt{t_p-\sqrt{\xi}x}}\cdot J\\
&=& \frac{2[1+o(1)]\exp(-c_p^2/2)
}{\pi\sqrt{t_p-\sqrt{\xi}c_p}}\frac{\exp\big(-\frac{1}{2}(t_p-\sqrt{\xi}x)(t_p-\sqrt{\xi}c_p)\big)}
{\sqrt{t_p-\sqrt{\xi}x}\cdot \big[\sqrt{\xi}(t_p-\sqrt{\xi}x)-2c_p\big]},
\end{eqnarray*}
proving \eqref{q(x)-approx-small-x}.  At this point, the proofs of (i), (ii) and (iii) are completed. \eop

\medskip






\medskip

\noindent{\it Proof of Lemma~\ref{estofq}.}  Observe that
\begin{equation}\label{III}
q(t_p)=Eq(\bar\bz_{1};t_p)=\frac{1}{\Phi(c_p)}\int^{c_p}_{-\infty}q(x;t_p)\phi(x)dx.
\end{equation}
For all $x\in \mathbb{R}$, we have
\[
x^2+\Big(\frac{4-\xi}{4}t_p^2-\frac{\xi^2}{4}x^2-\frac{(2-\xi)\sqrt{\xi}}{2}xt_p
\Big)=\frac{4-\xi^2}{4}\Big(x-\frac{\sqrt{\xi}t_p}{2+\xi}\Big)^2+\frac{2}{2+\xi}t_p^2
\]
and therefore
\begin{eqnarray}\label{intermediate}
&&\exp\Big(-\frac12\Big(\frac{4-\xi}{4}t_p^2-\frac{\xi^2}{4}x^2-\frac{(2-\xi)\sqrt{\xi}}{2}xt_p\Big)
\Big)\exp\Big(-\frac12 x^2\Big)\nonumber\\
&=&\exp\Big(-\frac{t_p^2}{2+\xi}\Big)\exp\Big(\frac{\xi^2-4}{8}\Big(x-\frac{\sqrt{\xi}t_p}{2+\xi}\Big)^2\Big).
\end{eqnarray}
Then it follows from \eqref{q(x)-approx} that
\begin{eqnarray}\label{q*phi}
&&q(x;t_p)\phi(x)\nonumber\\
&=&\frac{[1+o(1)]\exp\big(-\frac{t_p^2}{2+\xi}\big)}{\pi\sqrt{\big(t_p-\sqrt{\xi}x\big)\cdot\big[(1-0.5\xi)t_p+0.5\xi^{3/2}x\big]}}
\exp\Big(\frac{\xi^2-4}{8}\Big(x-\frac{\sqrt{\xi}t_p}{2+\xi}\Big)^2\Big) ~~~~~~~
\end{eqnarray}
uniformly over $x_p\le x\le c_p$, where $x_p$ is a sequence of positive numbers satisfying the
conditions in Lemma~\ref{lem2}(ii). With this at hand, we next begin to prove (a), (b) and (c) stated in the lemma.

{\it (a)}. In this case,  $\xi<2$.   In view of \eqref{III}, rewrite
\begin{equation}\label{II}
q(t_p)=\frac{1}{\Phi(c_p)}\int^{c_p}_{x_p}q(x;t_p)\phi(x)dx+
\frac{1}{\Phi(c_p)}\int^{x_p}_{-\infty}q(x;t_p)\phi(x)dx :=J_1+J_2.
\end{equation}
 Define $\sigma^2=4/(4-\xi^2)$.  Then it follows from \eqref{q*phi}
that
\begin{equation*}
q(x; t_p)\phi(x)
=\frac{4[1+o(1)]\exp\Big(-\frac{t_p^2}{2+\xi}\Big)}{\sqrt{2\pi(4-\xi^2)}}
 \frac{1}{\sqrt{2\pi}\sigma}\frac{\exp\Big(-\frac{1}{2\sigma^2}\Big(x-\frac{\sqrt{\xi}t_p}{2+\xi}\Big)^2\Big)}
{\sqrt{\big(t_p-\sqrt{\xi}x\big)\cdot\big[(1-0.5\xi)t_p+0.5\xi^{3/2}x\big]}}
\end{equation*}
uniformly over $x_p\le x\le c_p$. By changing the variables  $x=s+\frac{\sqrt{\xi}t_p}{2+\xi}$, we
have
\begin{eqnarray}
\big(t_p-\sqrt{\xi}x\big)\cdot\big[(1-0.5\xi)t_p+0.5\xi^{3/2}x\big]
&=&\Big(\frac{2t_p}{2+\xi}
-\sqrt{\xi}s\Big)\cdot\Big[\Big(1-\frac{\xi}{2+\xi}\Big)t_p+0.5\xi^{3/2}s\Big]\nonumber\\
&=&\Big(\frac{2t_p}{2+\xi}
-\sqrt{\xi}s\Big)\Big(\frac{2t_p}{2+\xi}+0.5\xi^{3/2}s\Big).\label{fromx2s}
\end{eqnarray}
Due to \eqref{tp} and \eqref{cp}, it is easy to check that, there exist constants $\delta_i>0, i=1, 2, 3, 4$ satisfying
\begin{eqnarray}\label{tau2}
\tau_1:= x_p-\frac{\sqrt{\xi}t_p}{2+\xi}\le -\delta_1 c_p~~~\mbox{and}~~~
\tau_2:= c_p-\frac{\sqrt{\xi}t_p}{2+\xi}\ge \delta_2 c_p
\end{eqnarray}
and
\begin{equation}\label{interval}
\frac{2t_p}{2+\xi} -\sqrt{\xi}s\ge\delta_3t_p~~~\mbox{and}~~~
\frac{2t_p}{2+\xi}+0.5\xi^{3/2}s\ge \delta_4t_p
\end{equation}
uniformly over $\tau_1\le s\le\tau_2$ as $p$ is sufficiently large. Now we set
\begin{equation}\label{ff}
f(s)=\frac1{\sqrt{\big(\frac{2t_p}{2+\xi}
-\sqrt{\xi}s\big)\big(\frac{2t_p}{2+\xi}+0.5\xi^{3/2}s\big)}}, ~~\tau_1\le s\le \tau_2.
\end{equation}
Then we have $f(0)=(2+\xi)(2t_p)^{-1}$. From the estimates in \eqref{interval}, we have $|f'(s)|\le d/t_p^2$ uniformly for
$\tau_1\le s\le\tau_2$ for some constant $d>0$. By using the mean-value theorem, we obtain that $|f(s)-f(0)|\le d|s|/t_p^2$
uniformly over $\tau_1\le s\le\tau_2$ for all large $p$.  Changing the variables
$x=s+\sqrt{\xi}t_p(2+\xi)^{-1}$ in \eqref{q*phi}, we obtain
\begin{eqnarray*}
J_1&=&[1+o(1)]\int^{c_p}_{x_p}q(x; t_p)\phi(x)dx\\
&=&\frac{4[1+o(1)]\exp\big(-\frac{t_p^2}{2+\xi}\big)}{\sqrt{2\pi(4-\xi^2)}}
 \int^{c_p}_{x_p}\frac{1}{\sqrt{2\pi}\sigma}\frac{\exp\Big(-\frac{1}{2\sigma^2}\Big(x-\frac{\sqrt{\xi}t_p}{2+\xi}\Big)^2\Big)dx}
{\sqrt{\big(t_p-\sqrt{\xi}x\big)\big[(1-0.5\xi)t_p+0.5\xi^{3/2}x\big]}}\\
&=&\frac{4[1+o(1)]\exp\big(-\frac{t_p^2}{2+\xi}\big)}{\sqrt{2\pi(4-\xi^2)}}
 \int^{\tau_2}_{\tau_1}\frac{1}{\sqrt{2\pi}\sigma}\frac{\exp\big(-\frac{1}{2\sigma^2}s^2\big)ds}
{\sqrt{\big(\frac{2t_p}{2+\xi}
-\sqrt{\xi}s\big)\big(\frac{2t_p}{2+\xi}+0.5\xi^{3/2}s\big)}}\\
&=&\frac{4[1+o(1)]\exp\big(-\frac{t_p^2}{2+\xi}\big)}{\sqrt{2\pi(4-\xi^2)}}\Big\{
 \int^{\tau_2}_{\tau_1}\frac{2+\xi}{2t_p}\cdot\frac{\exp\big(-\frac{1}{2\sigma^2}s^2\big)ds}
{\sqrt{2\pi}\sigma}\\
&&~~~+\int^{\tau_2}_{\tau_1}\big[]
f(s)-f(0)\big]\frac{\exp\Big(-\frac{1}{2\sigma^2}s^2\Big)ds}
{\sqrt{2\pi}\sigma}\Big\}.
\end{eqnarray*}
Use the estimate $|f(s)-f(0)|\le d|s|/t_p^2$
for $\tau_1\le s\le\tau_2$ to see that the above is identical to
\begin{eqnarray*}
&&\frac{4[1+o(1)]\exp\big(-\frac{t_p^2}{2+\xi}\big)}{\sqrt{2\pi(4-\xi^2)}}\Big\{
 \frac{2+\xi}{2t_p}\Big[\Phi\Big(\frac{\tau_2}{\sigma}\Big)-\Phi\Big(\frac{\tau_1}{\sigma}\Big)\Big]\\
&&  ~~~~~~~~~~~~~~~~~~~~~~~~~~~~~~~~~~~~~~~~~~~~~~~~~~~+O\left(\frac{1}{t_p^2}\right)\int^{\tau_2}_{\tau_1}|s|\frac{\exp\big(-\frac{1}{2\sigma^2}s^2\big)ds}
{\sqrt{2\pi}\sigma}\Big\}\\
&=&\frac{4[1+o(1)]\exp\Big(-\frac{t_p^2}{2+\xi}\Big)}{\sqrt{2\pi(4-\xi^2)}}\Big\{
 \frac{2+\xi}{2t_p}[1+o(1)]+O\left(\frac{1}{t_p^2}\right)\int^{\infty}_{-\infty}|s|\frac{\exp\big(-\frac{1}{2\sigma^2}s^2\big)ds}
{\sqrt{2\pi}\sigma}\Big\}\\
&=&[1+o(1)]\sqrt{\frac{2(2+\xi)}{\pi(2-\xi)}}\frac{\exp\big(-\frac{t_p^2}{2+\xi}\big)}{t_p}
\end{eqnarray*}
as $p\to\infty$. We have used the facts that $\lim\limits_{p\to\infty}\tau_1=-\infty$
and $\lim\limits_{p\to\infty}\tau_2=\infty$ from \eqref{tau2}.

To estimate $J_2$, we have from \eqref{q(x)-upper}, \eqref{intermediate} and the definition of $\sigma^2$ that
\begin{equation*}
q(x; t)\phi(x)\le
\frac{4[1+o(1)]\exp\big(-\frac{t_p^2}{2+\xi}\big)}{\sqrt{2\pi(4-\xi^2)}\big(t_p-\sqrt{\xi}c_p\big)}
\frac{1}{\sqrt{2\pi}\sigma}\exp\Big(-\frac{1}{2\sigma^2}\Big(x-\frac{\sqrt{\xi}t_p}{2+\xi}\Big)^2\Big)
\end{equation*}
uniformly over $x\le x_p$. Then we get
\begin{eqnarray*}
J_2&\le&
[1+o(1)]\int^{x_p}_{-\infty}q(x;t_p)\phi(x)dx\\
&=&O\Big(\frac{1}{t_p}\exp\Big(-\frac{t_p^2}{2+\xi}\Big)\Big)\int^{x_p}_{-\infty}
\frac{1}{\sqrt{2\pi}\sigma}\exp\Big(-\frac{1}{ 2\sigma^2}\Big(x-\frac{\sqrt{\xi}t_p}{2+\xi}\Big)^2\Big)dx\\
&=&O\Big(\frac{1}{t_p}\exp\Big(-\frac{t_p^2}{2+\xi}\Big)\Big)\int^{\tau_1}_{-\infty}
\frac{1}{\sqrt{2\pi}\sigma}\exp\Big(-\frac{s^2}{2 \sigma^2}\Big)ds\\
&=&O\Big(\frac{1}{t_p}\exp\Big(-\frac{t_p^2}{2+\xi}\Big)\Big)\Phi\big(\tau_1\sigma^{-1}\big)\\
&=&o\Big(\frac{1}{t_p}\exp\Big(-\frac{t_p^2}{2+\xi}\Big)\Big)
\end{eqnarray*}
since $\Phi(\tau_1/\sigma)\to 0$ by \eqref{tau2}. Therefore, \eqref{case1} is obtained by combining \eqref{II} and
estimates of $J_1$ and $J_2$ above.

{\it (b)}. In this case, $\xi=2$.
In view of \eqref{III}, rewrite
\begin{eqnarray*}
&&q(t_p) \nonumber\\
&=&\frac{1}{\Phi(c_p)}\Big\{\int^{c_p}_{x_p}q(x;t_p)\phi(x)dx+\int^{x_p}_{2^{-1/2}t_p-c_p}q(x;t_p)\phi(x)dx+
\int^{2^{-1/2}t_p-c_p}_{-\infty}q(x;t_p)\phi(x)dx\Big\}\nonumber\\
&=&:J_3+J_4+J_5,
\end{eqnarray*}
where $x_p$ and $b_p$ are two sequences of positive numbers satisfying the
restrictions in Lemma~\ref{lem2}(ii). By using \eqref{q*phi} and  changing the variables
$x=t_ps/\sqrt{2}$, we have
\begin{eqnarray*}
J_3
&=&\frac{[1+o(1)]\exp\big(-t_p^2/4\big)}{\pi}\int^{c_p}_{x_p}\frac{1}{\sqrt{(t_p-\sqrt{2}x\big)\sqrt{2}x}}dx\\
&=&\frac{[1+o(1)]\exp\big(-t_p^2/4\big)}{ \sqrt{2}\pi}\int^{\sqrt{2}c_p/t_p}_{\sqrt{2}x_p/t_p}s^{-1/2}(1-s)^{-1/2}ds \\
&=&
\frac{[1+o(1)]\exp\big(-t_p^2/4\big)}{\sqrt{2}\pi}\int^{\sqrt{2}c_p/t_p}_{1-\sqrt{2}c_p/t_p+\sqrt{2}b_p/t_p}s^{-1/2}(1-s)^{-1/2}ds\\
&=&\frac{[1+o(1)]\exp\big(-t_p^2/4\big)}{\sqrt{2}}\frac1{\pi}\int^{\sqrt{2}/r}_{1-\sqrt{2}/r}s^{-1/2}(1-s)^{-1/2}ds,
\end{eqnarray*}
where in the last step we use the fact
$0<\sqrt{2}b_p/t_p<b_p/c_p\to 0$ as $p\to\infty$. It is interesting to note
$\frac{1}{\pi}s^{-1/2}(1-s)^{1-/2}$, $0<s<1$, is the density function
of Beta($1/2$, $1/2$) distribution.

To show \eqref{case2}, it suffices to show that
$J_4=o(\exp(-t_p^2/4))$, and $J_5=o(\exp(-t_p^2/4))$. In fact, notice that, as $p$ is sufficiently large, $x_p=2^{-1/2}t_p-c_p+b_p<t_p$  and $t_p-\sqrt{2}c_p>\delta_5c_p$ by \eqref{tp} and \eqref{cp}, where $\delta_5$ is a constant free of $p$. It follows from \eqref{q(x)-upper} that
\begin{eqnarray*}
J_4&\le&
\frac{2[1+o(1)]\exp(-t_p^2/4)}{\sqrt{2\pi}\sqrt{t_p-\sqrt{2}c_p}}\int^{x_p}_{2^{-1/2}t_p-c_p}\frac{dx}{\sqrt{t_p-\sqrt{2}x}}\\
&\le&
\frac{2[1+o(1)]\exp(-t_p^2/4)}{\sqrt{2\pi}(t_p-\sqrt{2}c_p)}\Big[x_p-\Big(\frac{t_p}{\sqrt{2}}-c_p\Big)\Big]\\
&\le&
\frac{2[1+o(1)]\exp(-t_p^2/4)}{\sqrt{2\pi}\delta_5}\frac{b_p}{c_p}\\
&=&o\big(\exp(-t_p^2/4)\big).
\end{eqnarray*}
Review $\bar c_p(t_p,x)=c_p-t_p/\sqrt{2}+x$.
 It follows from  \eqref{q(x)-upper} that
\begin{eqnarray*}
J_5&\le&
\frac{2[1+o(1)]\exp(-t_p^2/4)}{\sqrt{2\pi}(t_p-\sqrt{2}c_p)}\int^{2^{-1/2}t_p-c_p}_{-\infty}\Phi(\bar
c_p(t_p,x))dx\\
&=&o\big(\exp(-t_p^2/4)\big)\int^{2^{-1/2}t_p-c_p}_{-\infty}\int^{c_p-2^{-1/2}t_p+x}_{-\infty}\phi(s)dsdx\\
&=&o\big(\exp(-t_p^2/4)\big)\int^{\infty}_{-\infty}|s|\phi(s)ds\\
&=&o\big(\exp(-t_p^2/4)\big)
\end{eqnarray*}
as $n\to\infty$. In the third step,  we choose to integrate out the
variable $x$ first for the double integral above. The completes the proof of \eqref{case2}.

{\it (c)}. In this case,  $\xi>2$. First, we write $q(t_p)$ as the sum of four integrals
\begin{eqnarray}\label{KKK}
q(t_p)
&=&\frac{1}{\Phi(c_p)}\left\{\int^{c_p}_{c_p-1}+\int^{c_p-1}_{\frac{\sqrt{\xi}t_p}{2+\xi}}
+\int^{\frac{\sqrt{\xi}t_p}{2+\xi}}_{\bar x_p}+\int^{\bar x_p}_{-\infty}\right\}q(x;t_p)\phi(x)dx\nonumber\\
&=&:K_1+K_2+K_3+K_4,
\end{eqnarray}
where $b_p$ in the definition of $\bar x_p$ satisfies the condition
in Lemma~\ref{lem2}(iii). In the following we will show that $K_1$ and $K_4$ are asymptotically equal;  $K_2$ and $K_3$  are negligible.

Recall that $\tau_2=c_p-\frac{\sqrt{\xi}t_p}{2+\xi}$ as defined in
\eqref{tau2}. We know
\begin{equation}\label{tau0}
\tau_3:=\frac{\sqrt{\xi}t_p}{2+\xi}-\bar
x_p=\frac{2}{\xi}\left(c_p-\frac{\sqrt{\xi}t_p}{2+\xi}\right)+b_p=\frac{2+o(1)}{\xi}\tau_2
\end{equation}
as $p$ is sufficiently large.
From \eqref{gongming} we have $0<\tau_3<\tau_2-1$ for all large $p$. Let $c>0$ be a fixed number.  One can easily verify that
\begin{equation*}
\int^z_0e^{cs^2}ds=\frac{1+o(1)}{2cz}e^{cz^2}
\end{equation*}
as $z\to\infty.$ This implies that
\begin{equation}\label{int2}
\int^{z-1}_0e^{cs^2}ds=
o\Big(\frac{e^{cz^2}}{z}\Big)~~\mbox{and}~~\int^z_{z-1}e^{cs^2}ds=\frac{[1+o(1)]e^{cz^2}}{2cz}
\end{equation}
as $z\to\infty.$ Now we will employ Lemma~\ref{lem2} to estimate $K_1$, $K_2$, $K_3$
and $K_4$ one by one. The variable change $x=s+\frac{\sqrt{\xi}t_p}{2+\xi}$ is applied to the three terms $K_1$, $K_2$ and
$K_3$.

To study $K_1$, set
\begin{equation}\label{meinv}
h(x)=\sqrt{\big(t_p-\sqrt{\xi}x\big)\big[(1-0.5\xi)t_p+0.5\xi^{3/2}x\big]},~~ x_p\le x\le c_p,
\end{equation}
where $x_p=\frac1{\sqrt{\xi}}\big(t_p-\frac{2c_p}{\sqrt{\xi}}\big)+b_p$ is as in Lemma \ref{lem2}(ii). Then from \eqref{fromx2s},
\[
h\Big(s+\frac{\sqrt{\xi}t_p}{2+\xi}\Big)=\frac{1}{f(s)},~~ \tau_1\leq s \leq \tau_2,
\]
where $f(s)$ is defined in \eqref{ff} for the proof of (a).
As in the proof for part (a), we have from the mean value theorem that
$|f(s)-f(\tau_2)|\le d/t_p^2$ for
$s\in [\tau_2-1, \tau_2]$, where $d$ is the same constant as that below \eqref{ff}. This implies
\[
f(s)=f(\tau_2)[1+o(1)]
\]
uniformly over $\tau_2-1\le s\le \tau_2$.
It follows from \eqref{ff} that
\[
f(\tau_2)=\frac1{h(c_p)}=\frac{1}{\sqrt{\big(t_p-\sqrt{\xi}c_p\big)\big[(1-0.5\xi)t_p+0.5\xi^{3/2}c_p\big]}}.
\]
Therefore, in virtue of  \eqref{q*phi},  \eqref{tau2} and
\eqref{int2}, we have
\begin{eqnarray}
K_1&=& \frac{[1+o(1)]\exp\big(-\frac{t_p^2}{2+\xi}\big)}{\pi}
\int^{c_p}_{c_p-1}\frac{\exp\Big(\frac{\xi^2-4}{8}\big(x-\frac{\sqrt{\xi}t_p}{2+\xi}\big)^2\Big)}{h(x)}
dx \nonumber\\
&=&\frac{[1+o(1)]\exp\big(-\frac{t_p^2}{2+\xi}\big)}{\pi}
\int^{\tau_2}_{\tau_2-1}\frac{
\exp\Big(\frac{\xi^2-4}{8}s^2\Big)}{h(
s+\frac{\sqrt{\xi}t_p}{2+\xi})}
ds \nonumber\\
&=&\frac{[1+o(1)]\exp\big(-\frac{t_p^2}{2+\xi}\big)}{\pi h(c_p)}
\int^{\tau_2}_{\tau_2-1}\exp\Big(\frac{\xi^2-4}{8}s^2\Big)ds \nonumber \\
&=&\frac{[1+o(1)]\exp\big(-\frac{t_p^2}{2+\xi}\big)}{\pi h(c_p)}
\frac{4}{\xi^2-4}\frac{\exp\big(\frac{\xi^2-4}{8}\tau_2^2
\big)}{\tau_2} \label{kuailede}\\
&=&\frac{4[1+o(1)]}{(\xi^2-4)\tau_2}q(c_p;t_p)\phi(c_p) \nonumber\\
&=&\frac{4[1+o(1)]q(c_p;t_p)\phi(c_p)}{(\xi^2-4)c_p-(\xi-2)\sqrt{\xi}t_p}, \nonumber
\end{eqnarray}
where we have used \eqref{q*phi} with
$x=c_p$ to obtain the last two expressions.

Now we work on $K_2$. By using the same argument, we have
\begin{eqnarray*}
K_2&=&\frac{[1+o(1)]\exp\big(-\frac{t_p^2}{2+\xi}\big) }{\pi}
\int^{\tau_2-1}_0\frac{\exp\Big(\frac{\xi^2-4}{8}s^2\Big)}{h\big(
s+\frac{\sqrt{\xi}t_p}{2+\xi}\big)}ds\\
&\le&\frac{[1+o(1)]\exp\big(-\frac{t_p^2}{2+\xi}\big)}{\pi h(c_p)} \int^{\tau_2-1}_0\exp\Big(\frac{\xi^2-4}{8}s^2\Big)ds\\
&=&\frac{[1+o(1)]\exp\big(-\frac{t_p^2}{2+\xi}\big)}{\pi h(c_p)}o\Big(\frac{\exp\big(\frac{\xi^2-4}{8}\tau_2^2\big)}{\tau_2}\Big)\\
&=&o(K_1).
\end{eqnarray*}

Now consider $K_3$. Note $\frac{\sqrt{\xi}t_p}{2+\xi} \leq c_p.$ By applying Lemma~\ref{lem2}(i) and \eqref{intermediate} and then employing changes of variables, we get
\begin{eqnarray*}
K_3 &\le&
\frac{O(1)\exp\big(-\frac{t_p^2}{2+\xi}\big)}{t_p-\sqrt{\xi}c_p}
\int^{\frac{\sqrt{\xi}t_p}{2+\xi}}_{\bar
x_p}\exp\Big(\frac{\xi^2-4}{8}\Big(x-\frac{\sqrt{\xi}t_p}{2+\xi}\Big)^2\Big)dx\\
&=&\frac{O(1)\exp\big(-\frac{t_p^2}{2+\xi}\big)}{t_p-\sqrt{\xi}c_p}
\int^{0}_{-\tau_3}\exp\Big(\frac{\xi^2-4}{8}s^2\Big)ds\\
&=&\frac{O(1)\exp\Big(-\frac{t_p^2}{2+\xi}\Big)}{t_p-\sqrt{\xi}c_p}
\int^{\tau_3}_0\exp\Big(\frac{\xi^2-4}{8}s^2\Big)ds.
\end{eqnarray*}
Since $\tau_3 \leq \tau_2-1$ as aforementioned, the above is controlled by
\begin{eqnarray*}
&&\frac{O(1)\exp\big(-\frac{t_p^2}{2+\xi}\big)}{t_p-\sqrt{\xi}c_p}
\int^{\tau_2-1}_0\exp\Big(\frac{\xi^2-4}{8}s^2\Big)ds\\
&= &\frac{O(1)\exp\big(-\frac{t_p^2}{2+\xi}\big)}{t_p-\sqrt{\xi}c_p}
 o\Big(\frac{\exp\big(\frac{\xi^2-4}{8}\tau_2^2\big)}{\tau_2}\Big)\\
&= &
 \frac{h(c_p)}{t_p-\sqrt{\xi}c_p}\cdot o\big(K_1\big) \\
&=&o(K_1).
\end{eqnarray*}
Here we have used \eqref{tau0}, \eqref{int2}, \eqref{kuailede} and the fact that both
$h(c_p)$ and $t_p-\sqrt{\xi}c_p$ are of the same order $c_p$.


To estimate $K_4$, we will employ \eqref{q(x)-approx-small-x}.  In order
to more conveniently cite some equations established earlier, we
change the variable $x$ in the integrand of $K_4$ to $y$, that is,
\begin{eqnarray*}
K_4&=&\frac{1}{\Phi(c_p)}\int^{\bar x_p}_{-\infty}q(y;t_p)\phi(y)dy\\
&=&\frac{2(1+o(1)\exp(-\frac12c_p^2)}{\pi\sqrt{t_p-\sqrt{\xi}c_p}}
\int^{\bar
x_p}_{-\infty}\frac{\exp\Big(-\frac{1}{2}\big(t_p-\sqrt{\xi}c_p\big)\big(t_p-\sqrt{\xi}y\big)-\frac{1}{2}y^2\Big)}
{\sqrt{2\pi}\sqrt{t_p-\sqrt{\xi}y}~\big[\sqrt{\xi}(t_p-\sqrt{\xi}y)-2c_p\big]}dy.
\end{eqnarray*}
When $\xi>2$, we have from \eqref{tp} and \eqref{cp} that there exists some $\delta_6>0$ satisfying
\begin{equation*}
\tau_4:=\bar
x_p-\frac{\sqrt{\xi}}{2}\left(t_p-\sqrt{\xi}c_p\right)=\frac{\xi^2-4}{2\xi}\left(c_p-\frac{\sqrt{\xi}t_p}{2+\xi}\right)-b_p\ge
\delta_6 c_p
\end{equation*}
for all large $p$.
Now we apply the identity \eqref{wecan} with $x=c_p$ to have
\begin{eqnarray*}
K_4&=&\frac{2[1+o(1)]}{\pi\sqrt{t_p-\sqrt{\xi}c_p}}\exp\Big(-\frac12\Big(\frac{4-\xi}{4}t_p^2+\frac{4-\xi^2}{4}c_p^2-\frac{(2-\xi)\sqrt{\xi}}{2}c_pt_p\Big)
\Big)\\
&&\cdot\int^{\bar
x_p}_{-\infty}\frac{\exp\Big(-\frac12\big(y-\frac12\sqrt{\xi}\big(t_p-\sqrt{\xi}c_p\big)\big)^2\Big)}
{\sqrt{2\pi}\sqrt{t_p-\sqrt{\xi}y}~\big[\sqrt{\xi}\big(t_p-\sqrt{\xi}y\big)-2c_p\big]}
dy\\
&=&\frac{2[1+o(1)]}{\pi\sqrt{t_p-\sqrt{\xi}c_p}}\exp\Big(-\frac12\Big(\frac{4-\xi}{4}t_p^2+\frac{4-\xi^2}{4}c_p^2-\frac{(2-\xi)\sqrt{\xi}}{2}c_pt_p\Big)
\Big)\\
&&\cdot\int^{\tau_4}_{-\infty}\frac{\exp\big(-\frac12s^2\big)ds}
{\sqrt{2\pi}\sqrt{t_p-\sqrt{\xi}\big(\frac{\sqrt{\xi}}{2}\big(t_p-\sqrt{\xi}c_p\big)+s\big)}
~\Big[\sqrt{\xi}\big(t_p-\sqrt{\xi}\big(\frac{\sqrt{\xi}}{2}\big(t_p-\sqrt{\xi}c_p\big)+s\big)\big)-2c_p\Big]}
\\
&=&\frac{2[1+o(1)]}{\pi\sqrt{t_p-\sqrt{\xi}c_p}}\exp\Big(-\frac12\Big(\frac{4-\xi}{4}t_p^2+\frac{4-\xi^2}{4}c_p^2-\frac{(2-\xi)\sqrt{\xi}}{2}c_pt_p\Big)
\Big)\\
&&\cdot\int^{\tau_4}_{-\infty}\frac{\exp\big(-\frac12s^2\big)ds}
{\sqrt{2\pi}\sqrt{\big(1-0.5\xi\big)t_p+0.5\xi^{3/2}c_p-\sqrt{\xi}s}~
\Big[\frac{\xi^2-4}{2}\big(c_p-\frac{\sqrt{\xi}t_p}{2+\xi}\big)-\xi
s\Big]},
\end{eqnarray*}
where
in the second step, we
have changed the variables
$y=\frac{\sqrt{\xi}}{2}(t_p-\sqrt{\xi}c_p)+s$. By following the same
procedure for the proof of \eqref{proofofI2}, we see that the above
integral is equal to
\[
\frac{1+o(1)} {\sqrt{(1-0.5\xi)t_p+0.5\xi^{3/2}c_p}\cdot
\frac{\xi^2-4}{2}\big(c_p-\frac{\sqrt{\xi}t_p}{2+\xi}\big)}.
\]
The details are omitted here.   Then we conclude that
\[
K_4=\frac{4[1+o(1)]\exp\Big(-\frac12\Big(\frac{4-\xi}{4}t_p^2+\frac{(\xi-2)\sqrt{\xi}}{2}c_pt_p+\frac{4-\xi^2}{4}c_p^2\Big)\Big)
}{\pi(\xi^2-4)\sqrt{t_p-\sqrt{\xi}c_p}\sqrt{(1-0.5\xi)t_p+0.5\xi^{3/2}c_p\big)}
\cdot\big(c_p-\frac{\sqrt{\xi}t_p}{2+\xi}\big)},
\]
which is equal to
\[
\frac{4[1+o(1)]q(c_p;t_p)\phi(c_p)}{(\xi^2-4)c_p-(\xi-2)\sqrt{\xi}t_p}
\]
from \eqref{intermediate} and  \eqref{q*phi}. In summary, we have shown that $K_1$ and $K_4$ are asymptotically equal, and  $K_2$
and $K_3$ are asymptotically negligible. This and \eqref{KKK} conclude
\eqref{case3}.

To show \eqref{j-moments},  for any $j\ge 2$ we write
\begin{equation*}
E\big[q(\bar\bz_{1};t_p)^j\big]=\frac1{\Phi(c_p)}\left\{\int^{c_p}_{c_p-1}+\int^{c_p-1}_{-\infty}\right\}q(x;t_p)^j\phi(x)dx:=K_5+K_6.
\end{equation*}
We first estimate $K_6$ before studying $K_5.$

Under conditions \eqref{tp} and \eqref{cp}, we have
$t_p<(\sqrt{\xi}+\frac{2}{\sqrt{\xi}})c_p$ for all large $p$. Then
it follows from \eqref{q(x)-approx} that as $p\to\infty$,
\begin{eqnarray*}
\frac{q(c_p-1;t_p)}{q(c_p;t_p)}&=&[1+o(1)]e^{\xi^2/8}\exp\Big(\frac14(\xi-2)\sqrt{\xi}t_p-\frac14\xi^2c_p\Big)\\
&\le&[1+o(1)]e^{\xi^2/8}\exp\Big(\frac{1}4(\xi-2)\sqrt{\xi}\Big(\sqrt{\xi}+\frac{2}{\sqrt{\xi}}\Big)c_p-\frac{1}4\xi^2c_p\Big)\\
&\le&[1+o(1)]e^{\xi^2/8}e^{-c_p}\\
&\to& 0.
\end{eqnarray*}
 This implies
\[
K_6\le
\frac{q(c_p-1;t_p)^{j-1}}{\Phi(c_p)}\int^{c_p-1}_{-\infty}q(x;t_p)\phi(x)dx=o\Big(q(c_p;t_p)^{j-1}q(t_p)\Big).
\]
Define $\tau_5=c_p
-\frac{j(\xi-2)\sqrt{\xi}}{j\xi^2-4}t_p$. Since $j\ge 2$, we have from \eqref{tau2} that as $p\to\infty$,
\[
\tau_5\ge c_p
-\frac{(\xi-2)\sqrt{\xi}}{\xi^2-4}t_p=c_p
-\frac{\sqrt{\xi}}{\xi+2}t_p=\tau_2\to\infty.
\]
By following the same procedure for estimating $K_1$,  we have from \eqref{q(x)-approx} and \eqref{meinv} that
\begin{eqnarray*}
K_5 &=& \frac{2^j[1+o(1)]}{(\sqrt{2\pi})^{j+1}}
\int^{c_p}_{c_p-1}\frac{\exp\Big(-\frac{j}2\big(\frac{4-\xi}{4}t_p^2-\frac{\xi^2}{4}x^2-\frac{(2-\xi)\sqrt{\xi}}{2}xt_p\big)-\frac12
x^2\Big)}
{h(x)^j}dx\\
&=& \frac{2^j[1+o(1)]}{(\sqrt{2\pi})^{j+1}h(c_p)^j}
\int^{c_p}_{c_p-1}\exp\Big(\frac{j\xi^2-4}{8}x^2-\frac{j(4-\xi)}{8}t_p^2
-\frac{j(\xi-2)\sqrt{\xi}}{4}xt_p\Big)
dx\\
&=&
\frac{2^j[1+o(1)]\exp\Big(-\frac{t_p^2}{8}\Big(j(4-\xi) + \frac{j^2(\xi-2)^2\xi}{j\xi^2-4}
\Big)\Big)}{(\sqrt{2\pi})^{j+1}h(c_p)^j}\\
 &&~~\cdot\int^{c_p}_{c_p-1}\exp\Big(\frac{j\xi^2-4}{8}\Big(x
-\frac{j(\xi-2)\sqrt{\xi}}{j\xi^2-4}t_p\Big)^2\Big)dx.
\end{eqnarray*}
With a change of variable, the above is equal to
\begin{eqnarray*}
&& \frac{2^j[1+o(1)]\exp\Big(-\frac{t_p^2}{8}\Big(j(4-\xi) +\frac{j^2(\xi-2)^2\xi}{j\xi^2-4}
\Big)\Big)}{(\sqrt{2\pi})^{j+1}h(c_p)^j}\int^{\tau_5}_{\tau_5-1}\exp\Big(\frac{j\xi^2-4}{8}s^2\Big)ds\\
&=&
\frac{2^j[1+o(1)]\exp\Big(-\frac{t_p^2}{8}\Big(j(4-\xi) +\frac{j^2(\xi-2)^2\xi}{j\xi^2-4}
\Big)\Big)}{(\sqrt{2\pi})^{j+1}h(c_p)^j}\cdot\frac{4\exp\Big(\frac{j\xi^2-4}{8}\tau_5^2\Big)}{(j\xi^2-4)\tau_5}.
\end{eqnarray*}
One can verify via the definition of $\tau_5$ that
\begin{eqnarray*}
&&-\frac{t_p^2}{8}\Big(j(4-\xi)+\frac{j^2(\xi-2)^2\xi}{j\xi^2-4}\Big) +\frac{j\xi^2-4}{8}\tau_5^2\\
&=& -\frac12\Big(\frac{4-\xi}{4}t_p^2-\frac{\xi^2}{4}c_p^2-\frac{(2-\xi)\sqrt{\xi}}{2}c_pt_p\Big)\cdot j -\frac{1}{2}c_p^2.
\end{eqnarray*}
By comparing the above assertions with \eqref{q(x)-approx}, we know
\begin{eqnarray*}
K_5&=&\frac{4[1+o(1)]}{{(j\xi^2-4)\tau_5}}q(c_p;t_p)^j\phi(c_p)\\
&=&\frac{4[1+o(1)]}{(j\xi^2-4)c_p-j(\xi-2)\sqrt{\xi}t_p}q(c_p;t_p)^j\phi(c_p)
\end{eqnarray*}
as $p\to\infty$.  Then \eqref{j-moments} is obtained from the
estimates of $K_5$ and $K_6$ above.  This completes the proof. \eop

\medskip

\noindent{\it Proof of Lemma~\ref{chores1}.}
(i). Recall $\alpha_p=\sqrt{2\ln p}$ and \eqref{standard}. We have
\begin{equation*}
\beta_p=\alpha_p-\alpha_p^{-1}\ln\big(\sqrt{2\pi}\alpha_p\big)~~\mbox{
and}~~~c_p=c_p(y)=\beta_p+\alpha_p^{-1}y, ~~~y\in \mathbb{R}.
\end{equation*}
Then $c_p\sim \beta_p \sim \alpha_p=\sqrt{2\ln p}$ as $p\to\infty$ for any $y\in \mathbb{R}.$ Here we use the notation $f(p) \sim g(p)$ for two functions $f(p)$ and $g(p)$ if $f(p)/g(p) \to 1$ as $p\to \infty$. It follows that
\begin{eqnarray*}
    \frac{\phi(c_p)}{c_p}\sim \frac{1}{\sqrt{2\pi}\cdot\sqrt{2\ln p}}\cdot e^{-c_p^2/2}.
\end{eqnarray*}
Now $c_p^2=\beta_p^2+2\beta_p\alpha_p^{-1}y + o(1)=2\ln p-2\ln\big(\sqrt{2\pi}\cdot \sqrt{2\ln p}\big) + 2y +o(1)$. Therefore,
\begin{eqnarray*}
    \frac{\phi(c_p)}{c_p}=[1+o(1)] \frac{1}{2\sqrt{\pi\ln p}}\cdot e^{-\ln p+\ln\big(2\sqrt{\pi\ln p}\big) -y +o(1)}=\big(1+o(1)\big)\frac{e^{-y}}{p}.
\end{eqnarray*}
We obtain (i).

(ii). By definition,
\begin{equation*}
A_{p,\xi}=\frac{2}{\sqrt{2+\xi}}\alpha_p\ \ \mbox{and}\ \ B_{p,\xi}=\sqrt{2+\xi}\alpha_p-\frac{\sqrt{2+\xi}}{2\alpha_p}\ln\big(\sqrt{2\pi(2-\xi)}\alpha_p\big)
\end{equation*}
since $\xi \in (0, 2)$. Obviously, $t_p=B_{p,\xi}+A_{p,\xi}^{-1}z= \sqrt{2+\xi}\alpha_p[1+o(1)]$. For $\xi \in (0, 2)$, observe
\begin{eqnarray*}
t_p^2
&=&B_{p,\xi}^2+A_{p,\xi}^{-2}z^2 +2B_{p,\xi}A_{p,\xi}^{-1}z \nonumber\\
&=& (2+\xi)\alpha_p^2-(2+\xi)\ln\big(\sqrt{2\pi(2-\xi)}\alpha_p\big)
+ (2+\xi)z +o(1).
\end{eqnarray*}
We thus have from (a) of Lemma \ref{estofq} that
\begin{eqnarray*}
q(t_p)
&= & [1+o(1)]\sqrt{\frac{2(2+\xi)}{\pi(2-\xi)}}\cdot\frac{\exp\left(-\frac{t_p^2}{2+\xi}\right)}{t_p}\\
& = & [1+o(1)]\sqrt{\frac{2(2+\xi)}{\pi(2-\xi)}}\cdot\frac{1}{\sqrt{2(2+\xi)\ln p}}\cdot \exp\left(-2\ln p+\ln\left(\sqrt{4\pi(2-\xi)\ln p}\right)-z\right)\\
& = & [1+o(1)]\frac{2e^{-z}}{p^2}.
\end{eqnarray*}
The desired conclusion follows for $\xi \in (0, 2)$.

For $\xi=2$, trivially $r=\lim_{p\to\infty}t_p/\sqrt{2\ln p}=2$. According to Lemma \ref{estofq},
\begin{equation*}
 q(t_p)=\frac{[1+o(1)]\exp\big(-t_p^2/4\big)}{\sqrt{2}\pi}\int^{1/\sqrt{2}}_{1-(1/\sqrt{2})}s^{-1/2}(1-s)^{-1/2}ds.
 \end{equation*}
Now we evaluate this integral. In fact, set $\beta=\frac1\pi\int^{1/\sqrt{2}}_{1-(1/\sqrt{2})}s^{-1/2}(1-s)^{-1/2}ds$. We claim $\beta=(2/\pi)\arcsin(\sqrt{2}-1)$. In fact, set $t=\sqrt{s}$. Then
\begin{eqnarray*}
\pi \beta
=\int^{2^{-1/4}}_{\sqrt{1-2^{-1/2}}}\frac{1}{t\sqrt{1-t^2}}\cdot (2t)\,dt =2\arcsin (t)\Big|^{2^{-1/4}}_{\sqrt{1-2^{-1/2}}}.
\end{eqnarray*}
Therefore
\begin{eqnarray*}
\frac{\pi}{2}\beta= \arcsin \left(2^{-1/4}\right)-\arcsin \left(\sqrt{1-2^{-1/2}}\right).
\end{eqnarray*}
Here we understand $a:=\arcsin \left(2^{-1/4}\right)$ and $b:=\arcsin \left(\sqrt{1-2^{-1/2}}\right)$ are both in $(0, \pi/2).$ Since $a>b$, we know $a-b\in (0, \pi/2)$ and hence $\beta \in (0, 1)$. Use $\sin a=2^{-1/4}$ and $\sin b=\sqrt{1-2^{-1/2}}$ to see $\cos a=\sqrt{1-2^{-1/2}}$ and $\cos b=2^{-1/4}$, respectively. Therefore,
\begin{eqnarray*}
 \sin\left(\frac{\pi}{2}\beta\right) &=& \sin (a-b)\\
 &=&\sin(a)\cos (b)-\cos (a)\sin (b)\\
 &=& 2^{-1/4}\cdot 2^{-1/4} - \sqrt{1-2^{-1/2}}\cdot \sqrt{1-2^{-1/2}}\\
 &=& \sqrt{2}-1.
\end{eqnarray*}
This yields that $\beta=(2/\pi)\arcsin(\sqrt{2}-1).$ Therefore,
\begin{equation}\label{Wukong}
q(t_p)=[1+o(1)]\exp\big(-t_p^2/4\big) \cdot \frac{\sqrt{2}}{\pi}\arcsin(\sqrt{2}-1).
\end{equation}
By definition,  $A_{p,\xi}=\alpha_p$ and
\begin{equation*}
B_{p,\xi}=2\alpha_p-\frac{1}{\alpha_p}\ln\frac{\sqrt{2}\pi}{\arcsin(\sqrt{2}-1)}.
\end{equation*}
 Then
\begin{eqnarray*}
t_p^2
&=&B_{p,\xi}^2+A_{p,\xi}^{-2}z^2 +2B_{p,\xi}A_{p,\xi}^{-1}z \nonumber\\
&=& 4\alpha_p^2-4\ln\frac{\sqrt{2}\pi}{\arcsin(\sqrt{2}-1)}
+ 4z +o(1).
\end{eqnarray*}
This implies that, at $\xi=2$,
\[
\frac{1}{4}t_p^2=2\ln p-\ln\frac{\sqrt{2}\pi}{\arcsin(\sqrt{2}-1)}+z+o(1),
\]
which together with \eqref{Wukong} concludes
\begin{equation*}
q(t_p)=[1+o(1)]\frac{2e^{-z}}{p^2}.
\end{equation*}
In summary, we complete the proof of (ii).

(iii) We have by plugging $x=c_p$ in \eqref{q*phi} that
\begin{equation}\label{shadondong}
q(c_p;t_p)\phi(c_p)
=\frac{[1+o(1)]\exp\big(-\frac{t_p^2}{2+\xi}\big)}{\pi\sqrt{\big(t_p-\sqrt{\xi}c_p\big)\big[(1-0.5\xi)t_p+0.5\xi^{3/2}c_p\big]}}
\exp\left(\frac{\xi^2-4}{8}\left(c_p-\frac{\sqrt{\xi}t_p}{2+\xi}\right)^2\right).
\end{equation}
Now we focus on $\xi\in (0,2]$.  Under the settings in Lemma~\ref{chores1}, we have
\[\lim_{p\to\infty}\frac{t_p}{c_p}=\lim_{p\to\infty}\frac{B_{p,\xi}}{\alpha_p}=\sqrt{2+\xi},
\]
that is, \eqref{gongming} holds with $r=\sqrt{2+\xi}$. Then it follows that
\[
\lim_{p\to\infty}\frac{t_p-\sqrt{\xi}c_p}{t_p}=1-\frac{\sqrt{\xi}}{\sqrt{2+\xi}}\in (0,1),
\]
\[
\lim_{p\to\infty}\frac{(1-0.5\xi)t_p+0.5\xi^{3/2}c_p}{t_p}=1-0.5\xi
+\frac{0.5\xi^{3/2}}{\sqrt{2+\xi}}\in (0,\infty)
\]
and
\[
\lim_{p\to\infty}(c_p-\frac{\sqrt{\xi}t_p}{2+\xi})=\infty.
\]
These, coupled with \eqref{shadondong}, imply
\[
q(c_p;t_p)\phi(c_p)
=O\Big(\frac{1}{t_p}\exp\big(-\frac{t_p^2}{2+\xi}\big)\Big)
\exp\left(\frac{\xi^2-4}{8}\left(c_p-\frac{\sqrt{\xi}t_p}{2+\xi}\right)^2\right)
\]
and
\[
\exp\left(\frac{\xi^2-4}{8}\left(c_p-\frac{\sqrt{\xi}t_p}{2+\xi}\right)^2\right)=\begin{cases}
    o(1)& \mbox{ if } \xi\in (0,2);  \\
    1 & \mbox{ if } \xi=2.
\end{cases}
\]
According to parts (a) and (b) of Lemma ~\ref{estofq},
 $q(t_p)$ is of the order $\frac{1}{t_p}\exp\big(-\frac{t_p^2}{2+\xi}\big)$ if $\xi\in (0,2)$, and of the order
$\exp\big(-\frac{t_p^2}{2+\xi}\big)$ if $\xi=2$. Therefore, we conclude that
$q(c_p;t_p)\phi(c_p)=o(q(t_p))$ whenever $\xi\in (0, 2]$. This confirms (iii). \eop

\medskip

\noindent{\it Proof of Lemma~\ref{chores2}.}  (a). By definition,
\begin{eqnarray}
A_{p,\xi} &=& \frac{2+\sqrt{\xi}}{\xi+\sqrt{\xi}}\alpha_p=\eta\alpha_p;\label{kexiao1}\\
B_{p,\xi} &=& \frac{2+2\sqrt{\xi}+\xi}{2+\sqrt{\xi}}\beta_p-\alpha_p^{-1}\ln\frac{\sqrt{1+\sqrt{\xi}}}{2+\sqrt{\xi}}\nonumber\\
 &=& \frac{\eta+1}{\eta}\beta_p - \alpha_p^{-1}\ln\frac{\sqrt{1+\sqrt{\xi}}}{2+\sqrt{\xi}}.\label{kexiao2}
\end{eqnarray}
 Hence, $t_p=B_{p,\xi}+A_{p,\xi}^{-1}z= \eta^{-1}(\eta+1)\beta_p[1+o(1)]$. Consequently, use $c_p \sim \beta_p \sim \alpha_p=\sqrt{2\ln p}$ to see
\begin{eqnarray*}
(j\xi^2-4)c_p-j(\xi-2)\sqrt{\xi}t_p
&=&c_p\cdot \left[j\xi^2-4-j(\xi-2)\sqrt{\xi}\frac{2+2\sqrt{\xi}+\xi}{2+\sqrt{\xi}}\right][1+o(1)].
\end{eqnarray*}
Now
\begin{eqnarray*}
&&j\xi^2-4-j(\xi-2)\sqrt{\xi}\frac{2+2\sqrt{\xi}+\xi}{2+\sqrt{\xi}}\\
&=& j\sqrt{\xi}\left[\xi^{3/2}-(\xi-2)\frac{2+2\sqrt{\xi}+\xi}{2+\sqrt{\xi}}\right]-4\\
&=& (4j)\frac{\sqrt{\xi}(\sqrt{\xi}+1)}{2+\sqrt{\xi}} -4\\
&= & 4\frac{j-\eta}{\eta}.
\end{eqnarray*}
Combining the two assertions, we get
\begin{eqnarray*}
    (j\xi^2-4)c_p-j(\xi-2)\sqrt{\xi}t_p=4\eta^{-1}(j-\eta)c_p[1+o(1)].
\end{eqnarray*}

(b). First, note that
\begin{eqnarray*}
 \frac1{\sqrt{\xi}}\left(t_p-\frac{2c_p}{\sqrt{\xi}}\right)
 &=&  \frac1{\sqrt{\xi}}\left[\frac{\eta+1}{\eta}-\frac{2}{\sqrt{\xi}}+o(1)\right]c_p\\
 &=& \left[1-\frac{4}{\xi(2+\sqrt{\xi})} +o(1)\right]c_p.
\end{eqnarray*}
Evidently, the constant $1-4/[\xi(2+\sqrt{\xi})] \in (0, 1)$ for $\xi >2.$ Consequently,
\begin{eqnarray*}
x_p:=\frac{1}{\ln\ln (p+3)}\left[c_p-  \frac1{\sqrt{\xi}}\left(t_p-\frac{2c_p}{\sqrt{\xi}} \right)\right]
\end{eqnarray*}
is a positive sequence as $p$ is sufficiently large and $x_p/c_p\to 0$. We then apply \eqref{q(x)-approx} by choosing $x=c_p$ to have
\begin{equation}\label{shuoba}
q(c_p;t_p)=\frac{2[1+o(1)]\exp\left(-\frac12\Big(\frac{4-\xi}{4}t_p^2-\frac{\xi^2}{4}c_p^2-\frac{(2-\xi)\sqrt{\xi}}{2}c_pt_p\Big)
\right)}{\sqrt{2\pi}\sqrt{\big(t_p-\sqrt{\xi}c_p\big)\left[(1-0.5\xi)t_p+0.5\xi^{3/2}c_p\right]}}.
\end{equation}
By using the identity below \eqref{kexiao2} that $t_p/c_p=1+\eta^{-1} +o(1)$,  it is easy to verify that
\begin{eqnarray*}
(1-0.5\xi)t_p+0.5\xi^{3/2}c_p
&=&\frac{1}{2}c_p\cdot \left[(2-\xi)(1+\eta^{-1})+\xi^{3/2} + o(1)\right] \\
&=& \frac{1}{2}c_p\cdot \left[(2-\xi)\left(1+\frac{\xi + \sqrt{\xi}}{2+\sqrt{\xi}}\right)+\xi^{3/2} + o(1)\right]\\
&= & 2c_p\cdot\left[\frac{1+\sqrt{\xi}}{2+\sqrt{\xi}}+o(1)\right].
\end{eqnarray*}
Furthermore,
\begin{eqnarray*}
 t_p-\sqrt{\xi}c_p=c_p\cdot \left[1+\frac{\xi+\sqrt{\xi}}{2+\sqrt{\xi}}-\sqrt{\xi}+o(1)\right]  =(2c_p)\cdot\left[\frac{1}{2+\sqrt{\xi}} + o(1)\right].
\end{eqnarray*}
As a result,
\begin{equation}\label{yizhuang}
\sqrt{\big(t_p-\sqrt{\xi}c_p\big)\left[(1-0.5\xi)t_p+0.5\xi^{3/2}c_p\right]}
=(2c_p)\cdot \left[\frac{\sqrt{1+\sqrt{\xi}}}{2+\sqrt{\xi}}+o(1)\right].
\end{equation}
Now we turn to evaluate the exponent in \eqref{shuoba}. First,
\begin{eqnarray*}
&&    \frac{4-\xi}{4}t_p^2=\frac{4-\xi}{4}\left(B_{p,\xi}+A_{p,\xi}^{-1}z\right)^2=\frac{4-\xi}{4}\left[B_{p,\xi}^2 +2B_{p,\xi}A_{p,\xi}^{-1}z+o(1)\right];\\
&& B_{p,\xi}^2=\left(\frac{\eta+1}{\eta}\right)^2\beta_p^2-\frac{2(\eta+1)}{\eta}\cdot \ln\frac{\sqrt{1+\sqrt{\xi}}}{2+\sqrt{\xi}} +o(1);\\
&& 2B_{p,\xi}A_{p,\xi}^{-1}z=\frac{2(\eta+1)}{\eta^2}z +o(1)
\end{eqnarray*}
since $\beta_p/\alpha_p\to 1$. Thus,
\begin{equation*}
\frac{4-\xi}{4}t_p^2=\frac{4-\xi}{4}\left[\left(\frac{\eta+1}{\eta}\right)^2\beta_p^2-\frac{2(\eta+1)}{\eta}\cdot \ln\frac{\sqrt{1+\sqrt{\xi}}}{2+\sqrt{\xi}}+\frac{2(\eta+1)}{\eta^2}z + o(1)\right].
\end{equation*}
Recall \eqref{standard}. We have
\begin{equation*}
\beta_p=\alpha_p-\alpha_p^{-1}\ln\big(\sqrt{2\pi}\alpha_p\big)~~\mbox{
and}~~~c_p(y)=\beta_p+\alpha_p^{-1}y
\end{equation*}
for any $y\in \mathbb{R}$.
Then
\begin{eqnarray*}
\frac{\xi^2}{4}c_p^2 = \frac{\xi^2}{4}\left[\beta_p^2+ 2\beta_p\alpha_p^{-1}y +o(1)\right]=\frac{\xi^2}{4}\left[\beta_p^2+ 2y +o(1)\right].
\end{eqnarray*}
From \eqref{kexiao1} and \eqref{kexiao2},
\begin{eqnarray*}
\frac{(2-\xi)\sqrt{\xi}}{2}c_pt_p
&=& \frac{(2-\xi)\sqrt{\xi}}{2}
\left(\beta_p+\alpha_p^{-1}y\right)\cdot \left(\frac{\eta+1}{\eta}\beta_p - \alpha_p^{-1}\ln\frac{\sqrt{1+\sqrt{\xi}}}{2+\sqrt{\xi}} +\frac{z}{\eta \alpha_p}\right)\\
&=&\frac{(2-\xi)\sqrt{\xi}}{2}\left(\frac{\eta+1}{\eta}\beta_p^2-\ln\frac{\sqrt{1+\sqrt{\xi}}}{2+\sqrt{\xi}} +\frac{z}{\eta}+\frac{\eta+1}{\eta}y +o(1)\right).
\end{eqnarray*}
Join the above to see
\begin{eqnarray*}
&&\frac{4-\xi}{4}t_p^2-\frac{\xi^2}{4}c_p^2-\frac{(2-\xi)\sqrt{\xi}}{2}c_pt_p\\
&=&\left[\frac{4-\xi}{4}\left(\frac{\eta+1}{\eta}\right)^2-\frac{\xi^2}{4}-\frac{(2-\xi)\sqrt{\xi}}{2}\left(\frac{\eta+1}{\eta}\right)\right]\beta_p^2\\
&& +\left[-\frac{\xi^2}{2}  - \frac{(2-\xi)\sqrt{\xi}}{2} \cdot \frac{\eta+1}{\eta}\right]y +\left[\frac{4-\xi}{4}\cdot  \frac{2(\eta+1)}{\eta^2}-  \frac{(2-\xi)\sqrt{\xi}}{2}\cdot  \frac{1}{\eta}\right]z\\
&& +\left[- \frac{4-\xi}{4}\cdot \frac{2(\eta+1)}{\eta}  + \frac{(2-\xi)\sqrt{\xi}}{2}\right]\ln\frac{\sqrt{1+\sqrt{\xi}}}{2+\sqrt{\xi}}.
\end{eqnarray*}
Now we simplify the coefficients of the four terms on the right hand side.  We will use the identity $\eta=(2+\sqrt{\xi})(\xi+\sqrt{\xi})^{-1}$ repeatedly. First, the coefficient of $\beta_p^2$ is identical to
\begin{eqnarray*}
\frac{4-\xi}{4}\left(\frac{\eta+1}{\eta}\right)^2-\frac{\xi^2}{4}-\frac{(2-\xi)\sqrt{\xi}}{2}\left(\frac{\eta+1}{\eta}\right)= 1.
\end{eqnarray*}
The coefficient of $y$ is equal to
\begin{eqnarray*}
-\frac{\xi^2}{2}  - \frac{(2-\xi)\sqrt{\xi}}{2} \cdot \frac{\eta+1}{\eta}=-\frac{2}{\eta}.
\end{eqnarray*}
The coefficient of $z$ is equal to
\begin{eqnarray*}
\frac{4-\xi}{4}\cdot  \frac{2(\eta+1)}{\eta^2}-  \frac{(2-\xi)\sqrt{\xi}}{2}\cdot  \frac{1}{\eta}=\frac{2}{\eta}.
\end{eqnarray*}
The coefficient of $\ln\frac{\sqrt{1+\sqrt{\xi}}}{2+\sqrt{\xi}}$ is equal to
\begin{eqnarray*}
&&- \frac{4-\xi}{4}\cdot \frac{2(\eta+1)}{\eta}  + \frac{(2-\xi)\sqrt{\xi}}{2}\\
&=&-\eta\cdot \left[\frac{4-\xi}{4}\cdot  \frac{2(\eta+1)}{\eta^2}-  \frac{(2-\xi)\sqrt{\xi}}{2}\cdot  \frac{1}{\eta}\right]\\
&=& -2
\end{eqnarray*}
by using the expression of the coefficient of $z$. Now, we combine the above displays to see that
\begin{eqnarray*}
&& -\frac{1}{2}\left(\frac{4-\xi}{4}t_p^2-\frac{\xi^2}{4}c_p^2-\frac{(2-\xi)\sqrt{\xi}}{2}c_pt_p\right)\\
&=& -\frac{1}{2}\beta_p^2+\ln\frac{\sqrt{1+\sqrt{\xi}}}{2+\sqrt{\xi}} +\frac{1}{\eta} (y-z) +o(1)\\
&=& \ln\frac{1}{p} +\ln(\sqrt{2\pi}\alpha_p) +\ln\frac{\sqrt{1+\sqrt{\xi}}}{2+\sqrt{\xi}} +\frac{1}{\eta} (y-z) +o(1)
\end{eqnarray*}
Finally, by this identity, \eqref{shuoba} and \eqref{yizhuang} we arrive at
\begin{eqnarray*}
&& pq(c_p;t_p) \\
&=& p\cdot \frac{2}{\sqrt{2\pi}} \cdot \frac{1}{2c_p}\cdot \left(\frac{\sqrt{1+\sqrt{\xi}}}{2+\sqrt{\xi}}\right)^{-1}\cdot \left[\frac{1}{p}\cdot (\sqrt{2\pi}\alpha_p)\cdot \frac{\sqrt{1+\sqrt{\xi}}}{2+\sqrt{\xi}}\cdot e^{(y-z)/\eta}\right][1+o(1)] \\
&\to &  e^{(y-z)/\eta}
\end{eqnarray*}
as $p\to\infty$.

(c).
From \eqref{case3} we have that
\begin{eqnarray*}
q(t_p)=\frac{8[1+o(1)]q(c_p;t_p)\phi(c_p)}{(\xi^2-4)c_p-(\xi-2)\sqrt{\xi}t_p}=\frac{8[1+o(1)]q(c_p;t_p)}{(\xi^2-4)-(\xi-2)\sqrt{\xi}(t_p/c_p)}\cdot \frac{\phi(c_p)}{c_p}.
\end{eqnarray*}
Recall the identity below \eqref{kexiao2} that $t_p/c_p=1+\eta^{-1} +o(1)$. From Lemma \ref{chores1}(i) we have
$\phi(c_p)/c_p=\big(1+o(1)\big)e^{-y}/p$. Apply conclusion (b) to get
\begin{eqnarray*}
p^2q(t_p)\to \frac{8e^{(y-z)/\eta}}{(\xi^2-4)-(\xi-2)\sqrt{\xi}(1+\eta^{-1})}  \cdot e^{-y}.
\end{eqnarray*}
Now
\begin{eqnarray}
 (\xi^2-4)-(\xi-2)\sqrt{\xi}(1+\eta^{-1})
 &=&(\xi-2)\left(\xi+2-\sqrt{\xi}\frac{2+2\sqrt{\xi}+\xi}{2+\sqrt{\xi}}\right) \nonumber\\
 &=& \frac{4(\xi-2)}{2+\sqrt{\xi}}=\frac{4(1-\eta)}{\eta}.\label{zaochen}
\end{eqnarray}
In conclusion,
\begin{eqnarray*}
p^2q(t_p)\to \frac{2\eta e^{-y}\cdot e^{(y-z)/\eta}}{1-\eta}
\end{eqnarray*}
The desired conclusion follows since $\tau=e^{(y-z)/\eta}$.

(d).
First, $(j-\eta)^{-1}$ is well-defined for $j\ge 2$. This is because  $\eta=(2+\sqrt{\xi})(\xi+\sqrt{\xi})^{-1}<1$ due to the assumption $\xi>2.$ By \eqref{j-moments},
\begin{equation*}
E\left(q(\bar\bz_{1};t_p)^j\right)=\frac{4[1+o(1)]q(c_p;t_p)^{j}\phi(c_p)}{(j\xi^2-4)c_p-j(\xi-2)\sqrt{\xi}t_p}
+o\Big(q(c_p;t_p)^{j-1}q(t_p)\Big).
\end{equation*}
Second, in view of Lemma \ref{chores1}(i) we see
$\phi(c_p)/c_p=\big(1+o(1)\big)e^{-y}/p$. Then by applying  conclusion (b), we get
\begin{eqnarray*}
q(c_p;t_p)^{j-1}q(t_p) = O\left(\frac{1}{p^{j-1}}\cdot \frac{1}{p^2}\right)=O\left(\frac{1}{p^{j+1}}\right).
\end{eqnarray*}
The two displays imply
\begin{eqnarray*}
p^{j+1}E\big(q(\bar\bz_{1};t_p)^j\big)=   \frac{4[1+o(1)][pq(c_p;t_p)]^{j}}{(j\xi^2-4)-j(\xi-2)\sqrt{\xi}(t_p/c_p)}\cdot \frac{p\phi(c_p)}{c_p}
+o\Big(1\Big).
\end{eqnarray*}
Similar to the argument used in the proof of (c), we see
\begin{eqnarray*}
 p^{j+1}E\big(q(\bar\bz_{1};t_p)^j\big) \to   \frac{4\tau^{j}}{(j\xi^2-4)-j(\xi-2)\sqrt{\xi} (1+\eta^{-1})}  \cdot e^{-y}.
\end{eqnarray*}
Use \eqref{zaochen} to obtain
\begin{eqnarray*}
&& (j\xi^2-4)-j(\xi-2)\sqrt{\xi} (1+\eta^{-1}) \\
 &=& j\left[(\xi^2-4) -(\xi-2)\sqrt{\xi} (1+\eta^{-1})\right] +4(j-1)\\
 &=& j\frac{4(1-\eta)}{\eta}+4(j-1)\\
 &=& 4(j-\eta)\eta^{-1}.
\end{eqnarray*}
This concludes $p^jE\big(q(\bar\bz_{1};t_p)^j\big) \to \eta(j-\eta)^{-1}\tau^je^{-y}.$
\eop




\baselineskip 12pt

\end{document}